\documentclass[12pt,twoside]{amsart}
\usepackage{amssymb,amsmath,amsthm, amscd, enumerate, mathrsfs}
\usepackage{graphicx, hhline}
\usepackage[all]{xy}
\usepackage[usenames]{color}
\usepackage{hyperref}
\usepackage{fancyhdr}
\usepackage{bbm}
\usepackage[top=30truemm,bottom=30truemm,left=30truemm,right=30truemm]{geometry}
\usepackage{comment}
\usepackage{todonotes}

\hypersetup{colorlinks=true}

\title{Special K-stability and positivity of CM line bundles}
\author{Masafumi Hattori}
\date{\today}

\address{Masafumi Hattori \\ Department of 
Mathematics, Graduate School of Science, 
Kyoto University, Kyoto 606-8502, Japan}
\email{hattori.masafumi.47z@st.kyoto-u.ac.jp}

\newtheorem{thm}{Theorem}[section]

\newtheorem{lem}[thm]{Lemma}
\newtheorem{cor}[thm]{Corollary}
\newtheorem{prop}[thm]{Proposition}

\newtheorem{claim}{Claim}

\theoremstyle{definition}
\newtheorem{de}[thm]{Definition}
\newtheorem{ex}[thm]{Example}

\newtheorem*{ack}{Acknowledgments}

\newtheorem*{claim*}{Claim}
\begin{document}

\maketitle
\begin{abstract}
    We show that the CM line bundle on a proper family parametrizing specially K-stable varieties with maximal variation is ample.
    As an application, we show projectivity of any proper subspace of the coarse moduli space of uniformly adiabatically K-stable klt--trivial fibrations over curves constructed in \cite{HH}.
\end{abstract}

\section{Introduction}

\subsection{Positivity of CM line bundles}
K-stability was first introduced by Tian \cite{T} for Fano manifolds and reformulated by Donaldson \cite{Dn2} for other polarized varieties.
This notion is defined in a purely algebro-geometric way and is considered to be closely related to the existence of constant scalar curvature K\"{a}hler (for short cscK) metrics as the Yau-Tian-Donaldson conjecture predicts.
On the other hand, K-stability is closely related to birational geometry and moduli problems.
For log Fano pairs, K-stability is completely detected by the $\delta$-invariant, which was first introduced by Fujita--Odaka \cite{FO} and Blum--Jonsson \cite{BlJ} and the moduli scheme parametrizing K-polystable log Fano pairs is constructed (see \cite{Xu}).
On the other hand, Hashizume and the author \cite{HH} construct a moduli space parametrizing uniformly adiabatically K-stable klt--trivial fibrations over curves and show that this moduli is a kind of K-moduli spaces if we choose some polarizations.

Odaka (\cite[Conjecture 5.2]{O2}) conjectured the K-moduli conjecture ten years ago, which predicts that there would exist a quasi-projective moduli scheme parametrizing all K-polystable varieties with fixed some numerical data.
This conjecture is shown for the Fano case (cf.~\cite{XZ}, \cite{LXZ}, \cite{Xu}).
He further conjectured that we can choose a $\mathbb{Q}$-ample line bundle on the moduli space to be the CM line bundle, which was introduced by Paul--Tian \cite{PT}, \cite{PT2}. 
Indeed, the CM line bundle of a KSBA moduli space (see \cite{kollar-moduli} for details) is ample by \cite{PX}.
Recently, it was shown that a K-moduli space of $\mathbb{Q}$-Fano varieties with fixed dimension and volume is a projective scheme (\cite{CP}, \cite{P}, \cite{XZ}) with the CM line bundle ample.
On the other hand, Fujiki--Schumacher \cite{FS} showed that compact subspaces of moduli spaces parametrizing some cscK manifolds are projective by using the generalized Weil--Petersson metric, which is closely related to the CM line bundle.
Their result is enhanced by Dervan--Naumann \cite{DN} for any moduli parametrizing all cscK manifolds.
Recently, Ortu \cite{Ort} has constructed moduli spaces parametrizing manifolds with optimal symplectic connections (cf.~\cite{DS1}).
We note that her moduli spaces parametrize fibrations that do not necessarily have cscK metrics.
Ortu also shows that any compact subspace of her moduli spaces is projective.

\subsection{Main results}

In this paper, we settle the conjecture of Odaka on positivity of the CM line bundle for special K-stability.
\begin{thm}\label{thm--main--2}
    Let $\pi\colon(X,\Delta,L)\to S$ be a polarized log $\mathbb{Q}$-Gorenstein family with maximal variation (cf.~Definitions \ref{de--polarized--family} and \ref{de--maximal-variation}), where $S$ is projective and $(X_{\bar{s}},\Delta_{\bar{s}})$ is klt for any geometric point $\bar{s}\in S$.
    If $(X_{s},\Delta_{s},L_{s})$ is specially {\rm K}-stable for any closed point $s\in S$,
then the CM-line bundle $\lambda_{\mathrm{CM},\pi}$ is ample.
\end{thm}

Special K-stability was first introduced by \cite{CM}. 
We note that K-stability of $\mathbb{Q}$-Fano varieties is equivalent to special K-stability of them and special K-stability is compatible with the theory of filtrations as Ding stability of log Fano pairs (see \cite{Fjt} and \cite{Li3}).

We recall the following fact shown in \cite{HH}.
\begin{thm}[For details, see Theorem {\ref{quesmain}}]\label{thm--main--0}
We fix $d\in\mathbb{Z}_{>0}$, $u\in\mathbb{Q}_{>0}$, $v\in\mathbb{Q}_{>0}$.
Then there exists a separated Deligne-Mumford moduli stack $\mathscr{M}_{d,v,u,r}$ of finite type over $\mathbb{C}$ with a coarse moduli space $M_{d,v,u,r}$ parametrizing uniformly adiabatically {\rm K}-stable klt-trivial fibrations $f\colon(X,0,A)\to\mathbb{P}^1$ such that
\begin{enumerate}
\item $\mathrm{dim}\,X=d$,
\item for any general fiber $F$, $F\cdot A^{d-1}=v$, and
\item $K_X\sim_{\mathbb{Q}}-uf^*\mathcal{O}(1)$.
\end{enumerate}

Furthermore, there exists $w\in\mathbb{Q}_{>0}$ such that for any uniformly adiabatically {\rm K}-stable klt-trivial fibration $f\colon(X,0,A)\to \mathbb{P}^1$ as above, if $\mathrm{vol}(A)=w$, then $(X,A)$ is specially {\rm K}-stable.
\end{thm}

By Theorem \ref{thm--main--2}, we show that proper subspaces of the coarse moduli spaces constructed by \cite{HH} are all projective.
\begin{cor}\label{thm--main--1}
Any proper subspace $B$ of $M_{d,v,u,r}$ is projective.   
\end{cor}

\subsection{Outline of the proof}
We briefly explain the idea of the proof of Theorem \ref{thm--main--2} here.
Special K-stability of a polarized klt pair $(X,\Delta,L)$ consists of the following two properties.
\begin{enumerate}[(i)]
    \item $H:=\delta(X,\Delta,L)L+K_X+\Delta$ is ample, and\label{special--property--(i)}
    \item uniform $\mathrm{J}^H$-stability of $(X,L)$.
\end{enumerate}
Consider a polarized log family $\pi\colon (X,\Delta,L)\to S$ of relative dimension $n$ and an ample line bundle $H$ on $X$.
Recall that the log CM line bundle $\lambda_{\mathrm{CM},\pi}$ is defined as
\[
\pi_*\left(-\frac{n(K_{X_s}+\Delta_s)\cdot L_s^{n-1}}{L_s^n}L^{n+1}+(n+1)L^n\cdot(K_{X/S}+\Delta)\right)
\]
(see Definition \ref{de-cm-line-bundle}). We define the following variant of the CM line bundle
\[
\lambda_{\mathrm{J},\pi,H}:=\pi_*\left(-\frac{nH_s\cdot L_s^{n-1}}{L_s^n}L^{n+1}+(n+1)L^n\cdot H\right).
\]
We call this the $\mathrm{J}^H$-line bundle (cf.~Definition \ref{de--j--line}).
The following is a key observation in this paper.
\begin{prop}\label{prop--j--ample}
Let $\pi\colon(X,L)\to S$ be a polarized family with an ample $\mathbb{Q}$-line bundle $H$, where $S$ is projective and every geometric fiber of $\pi$ is normal.
   If $(X_s,L_s)$ is $\mathrm{J}^{H_s}$-semistable for any closed point $s\in S$,
   then $\lambda_{\mathrm{J},\pi,H}$ is ample.
\end{prop}

This phenomenon is first observed by Murakami \cite[Lemma 2.7]{Mur} for the case when $\pi$ is smooth and all fibers are uniformly J-stable.
In this paper, we show that the same phenomenon occurs also when $\pi$ is non-smooth but flat.
This is the first ingredient to show Theorem \ref{thm--main--2}.
On the other hand, we show that $\pi_*(H^{n+1})$ is ample by using a similar technique to obtain the ampleness of the CM line bundle of a family of K-stable log Fano pairs with maximal variation as follows.
\begin{thm}\label{thm--delta--ampleness}
    Let $\pi\colon(X,\Delta,L)\to S$ be a polarized log $\mathbb{Q}$-Gorenstein family of relative dimension $n$ with maximal variation, where $S$ is projective and $(X_{\bar{s}},\Delta_{\bar{s}})$ is klt for any geometric point $\bar{s}\in S$.
    Suppose that $\pi_*L^{n+1}\equiv0$.
    If there exists $\lambda\in\mathbb{Q}_{>0}$ such that $\lambda<\delta(X_s,\Delta_s,L_s)$ for any closed point $s\in S$ and $K_{X/S}+\Delta+\lambda L$ is $\pi$-ample, then the $\mathbb{Q}$-line bundle $\pi_*(K_{X/S}+\Delta+\lambda L)^{n+1}$ is ample.
\end{thm}
This is the second ingredient.
By these ingredients, we obtain Theorem \ref{thm--main--2}.

\begin{ack}
The author would like to thank Professor Yuji Odaka for careful reading his draft.
He would also like to thank Rei Murakami for pointing out some typos.
    This work is partially supported by JSPS KAKENHI 22J20059  
(Grant-in-Aid for JSPS Fellows DC1).
\end{ack}

\section{Preliminaries}\label{sec2}

We work over the field of complex numbers $\mathbb{C}$.

\subsection*{Notations and conventions}
\begin{enumerate}[(i)]
\item If we say that $X$ is a scheme, then we assume $X$ to be of finite type over $\mathbb{C}$.
If $X$ is further separated, irreducible and reduced, then we say that $X$ is a variety.
For any point $x\in X$, let $\kappa(x)$ denote the residue field of the local ring $\mathcal{O}_{X,x}$.
That is, if we set $\mathfrak{m}_x$ as the maximal ideal of $\mathcal{O}_{X,x}$, $\kappa(x):=\mathcal{O}_{X,x}/\mathfrak{m}_x$.
\item
Let $X$ be a scheme.
We denote 
\[
X(S):=\mathrm{Hom}(S,X)
\]
and call this the set of all $S$-valued points of $X$.
If $S=\mathrm{Spec}\,\Omega$, where $\Omega$ is an algebraically closed field, then we call elements of $X(S)$ geometric points of $X$.
If the image of $\mathrm{Spec}\,\Omega\to X$ is $x\in X$, we denote this by $\bar{x}\in X$. 
\item Let $X$ be a scheme of finite type over $\mathbb{C}$ and $U$ an open subset.
We say that $U$ is {\it big} if $\mathrm{codim}_{X}(X\setminus U)\ge2$.
\item Let $X$ be a scheme of finite type over $\mathbb{C}$ with a ($\mathbb{Q}$-)line bundle $H$.
Let $|H|$ (resp.~$|H|_{\mathbb{Q}}$) denote the set of all effective divisors linearly equivalent (resp.~effective $\mathbb{Q}$-divisors $\mathbb{Q}$-linearly equivalent) to $H$.
    \item Let $X$ be a scheme of finite type over $\mathbb{C}$.
    We say that a property $\mathcal{P}$ {\it holds for any very general closed point} $x\in X$ if there exist countably many closed subvarieties $V_{i}\subsetneq X$ such that $\mathcal{P}$ holds for any closed point $x\in X\setminus\bigcup_{i=1}^\infty V_{i}$.
    \item \label{Notations--(vi)} Let $f\colon X\to S$ be a proper morphism of normal varieties.
    Let $g\colon T\to S$ be a morphism from a normal variety.
     Then we set $X_T:=X\times_ST$.
     Let $h\colon X_T\to X$ be the canonical morphism.
     If $L$ is a $\mathbb{Q}$-line bundle on $X$, we set $L_T:=h^*L$.
     For any $s\in S$, we denote $X_s=X_T$ and $L_s=L_T$, where $T=\mathrm{Spec}\,\kappa(s)$.
Suppose that $f$ is flat, all geometric fibers are connected and normal, there exists a $\mathbb{Q}$-divisor $\Delta$ such that $K_{X/S}+\Delta$ is $\mathbb{Q}$-Cartier and $\Delta$ does not contain any fiber of $X$ over $S$.
Then, we set $\Delta_T$ as follows.
Take an open subset $U\subset X$ such that $\mathrm{codim}_{X_s}(X_s\setminus U)\ge2$ for any $s\in S$ and $f|_U$ is smooth.
Then we see that $\Delta|_U$ is $\mathbb{Q}$-Cartier and can consider $h|_{h^{-1}(U)}^*\Delta|_U$.
Let $\Delta_T$ be the closure of $h|_{h^{-1}(U)}^*\Delta|_U$.
We note that then
\[
K_{X_T/T}+\Delta_T=h^*(K_{X/S}+\Delta).
\]
\item Let $D$ be a $\mathbb{Q}$-Weil divisor on a projective normal variety $X$.
We say that $D$ is big if there exists an ample $\mathbb{Q}$-Cartier $\mathbb{Q}$-divisor $A$ and an effective $\mathbb{Q}$-Weil divisor $E$ such that $D=A+E$.
If $D+A'$ is big for any ample $\mathbb{Q}$-Cartier $\mathbb{Q}$-divisor $A'$, then we say that $D$ is {\it pseudo-effective}.
\item Let $X$ be a projective normal variety.
We say that $C$ is a {\it movable curve} of $X$ if there exist a projective birational morphism $\mu\colon X'\to X$ and very ample hypersurfaces $H_1,\ldots,H_{\text{dim}\,X-1}$ such that $C=H_1\cap\ldots\cap H_{\text{dim}\,X-1}$.
Let $L$ be a line bundle on $X$.
We denote $L\cdot C:=\mu^*L\cdot C$ for simplicity.
By \cite{BDPP}, $L\cdot C\ge0$ for any movable curve $C$ if and only if $L$ is pseudo-effective.
\item Let $f\colon X\to S$ be a proper morphism such that any geometric fiber is normal and connected.
Then we consider the following functor.
For any morphism of schemes $T\to S$, we attain the following set
\[
\mathfrak{Pic}_{X/S}(T):=\{L\,|\,\text{$L$ is a line bundle on }X_T\}/\sim_T
\]
where $L_1\sim_TL_2$ if and only if $L_1\otimes f_T^*B\sim L_2$ for some line bundle $B$ on $T$.
Then we have the relative Picard scheme $\mathrm{Pic}_{X/S}$, which represents the \'{e}tale sheafification of the above functor.
If $S=\mathrm{Spec}\,\mathbb{C}$, then we simply denote $\mathrm{Pic}(X):=\mathrm{Pic}_{X/S}$. 
Furthermore, $\mathrm{Pic}^0(X)$ denotes the identity component of $\mathrm{Pic}(X)$ and parametrizes all line bundles algebraically equivalent to $\mathcal{O}_X$.
See \cite[\S9]{FGA} for details.
Let $[L]$ denote an element of $\mathrm{Pic}(X)$ whose representative is a line bundle $L$.
\item Let $E$ be a locally free sheaf on a smooth projective curve $C$.
We say that $E$ is {\it nef} (resp.~{\it ample}) if $\mathcal{O}_{\mathbb{P}_C(E)}(1)$ is nef (resp.~ample).
$E$ is called {\it weakly positive} if for any $a\in\mathbb{Z}_{>0}$ and ample line bundle $A$, the stalk of  $\mathrm{Sym}^{ab}E\otimes \mathcal{O}_C(bA)$ at the generic point of $C$ is generated by $H^0(C,\mathrm{Sym}^{ab}E\otimes \mathcal{O}_C(bA))$ for some $b\in\mathbb{Z}_{>0}$.
\item Let $X$ be a proper normal variety and $\pi\colon Y\to X$ be a resolution of singularities of $X$.
$\mathrm{Alb}(Y)$ denotes the Albanese variety of $Y$.
Let $\beta\colon  Y\to\mathrm{Alb}(Y)$ be a canonical morphism.
Then it is well-known that there exists a canonical morphism $\alpha\colon X\to \mathrm{Alb}(Y)$ such that $\beta=\alpha\circ\pi$.
Thus, we denote $\mathrm{Alb}(Y)$ by $\mathrm{Alb}(X)$ and call this the {\it Albanese variety} of $X$.
We also call $\alpha$ an {\it Albanese morphism}. 
\item Let $X$ be a Noetherian scheme.
Let $\mathscr{F}^\bullet$ be a complex of $\mathcal{O}_X$-modules.
We say that $\mathscr{F}^\bullet$ is a {\it perfect complex} if  there exist a family of open subsets $\{U_i\}_{i=1}^r$ and bounded complexes $\mathscr{G}_i^\bullet$ of finite free $\mathcal{O}_{U_i}$-modules with a quasi isomorphism $\mathscr{G}_i^\bullet\to\mathscr{F}^\bullet|_{U_i}$ for each $i$.
On the other hand, let $E$ be a coherent sheaf on $X$.
$E$ is called an $\mathcal{O}_X$-module {\it of finite Tor-dimension} if for any $x\in X$, the stalk $E_x$ admits a resolution of finitely generated free $\mathcal{O}_{X,x}$-modules of finite length (\cite[p.~111]{GIT}). 
\end{enumerate}

\subsection{K-stability}
We first recall the fundamental concepts of birational geometry and K-stability.

\begin{de}
Let $X$ be a quasi-projective normal variety.
Suppose that $B$ is a $\mathbb{Q}$-divisor on $X$ such that $K_X+B$ is $\mathbb{Q}$-Cartier.
Then we call $(X,B)$ a {\it sublog pair}.
If $B$ is further effective, then we say that $(X,B)$ is a {\it log pair}.
For any prime divisor $E$ over $X$, choose a proper birational morphism $\pi\colon Y\to X$ such that $E$ is defined on $Y$.
Then we set the log discrepancy of $(X,B)$ with respect to $E$ as 
\[
A_{(X,B)}(E)=\mathrm{ord}_E(K_Y-\pi^*(K_X+B))+1.
\]
The above value is independent of the choice of $\pi$.
We say that $(X,B)$ is {\it subklt} (resp.~{\it sublc}) if $A_{(X,B)}(E)>0$ (resp.~$\ge0$) for any prime divisor $E$ over $X$.
If $E$ is further effective, then we say that $(X,B)$ is {\it klt} (resp.~{\it lc}).

For any coherent ideal $\mathfrak{a}$ of $X$ and rational number $r>0$, we consider the pair $(X,B+r\mathfrak{a})$.
We call $r\mathfrak{a}$ a $\mathbb{Q}$-ideal.
We define the {\it log discrepancy} of $(X,B+r\mathfrak{a})$ as follows.
Take $E$ a prime divisor over $X$ and $\pi\colon Y\to X$ such that $\pi$ is a log resolution of $(X,B)$ and $\mathfrak{a}$, i.e., $\mathrm{Ex}(\pi)+\pi^{-1}_*B+\pi^{-1}\mathfrak{a}$ is simple normal crossing and $Y$ is a smooth variety proper over $X$ such that $E$ is a prime divisor defined on $Y$.
Here, we note that $\pi^{-1}\mathfrak{a}$ is a Cartier divisor.
Then 
\[
A_{(X,B+r\mathfrak{a})}(E)=\mathrm{ord}_E(K_Y-\pi^*(K_X+B)-r\pi^{-1}\mathfrak{a})+1.
\]
We say that $(X,B+r\mathfrak{a})$ is {\it subklt} (resp.~{\it sublc}) if $A_{(X,B+r\mathfrak{a})}(E)>0$ (resp.~$\ge0$) for any prime divisor $E$ over $X$.
\end{de}

Let $(X,B)$ be a log pair.
We set \[
\mathrm{Aut}(X,B):=\{\sigma\in\text{Aut}\,X\,|\,\sigma_*B=B\}.
\]
It is well-known that $\mathrm{Aut}(X,B)$ is a group scheme and let $\mathrm{Aut}_0(X,B)$ be the identity component of $\mathrm{Aut}(X,B)$.
\begin{de}[Log canonical threshold]
Let $(X,B)$ be a log subpair and let $D$ be an effective $\mathbb{Q}$-Cartier $\mathbb{Q}$-divisor on $X$.
Take an arbitrary $\mathbb{Q}$-ideal $\mathfrak{a}$ on $X$.
Then we define the {\it log canonical threshold} for $(X,B)$ with respect to $D$ as 
\[
\mathrm{lct}(X,B;D)=\sup\{t\in\mathbb{Q}\,|\,(X,B+tD) \text{ is sublc}\}.
\]
On the other hand, we set the log canonical threshold of $(X,B)$ with respect to $\mathfrak{a}$ as
\[
\mathrm{lct}(X,B;\mathfrak{a})=\sup\{t\in\mathbb{Q}\,|\,(X,B+t\mathfrak{a}) \text{ is sublc}\}.
\]
\end{de}

\begin{de}
We say that $\mathfrak{a}_{\bullet}$ is a graded sequence of non-zero ideals if there exists a sequence of ideals $\{\mathfrak{a}_m\}_{m\in\mathbb{Z}_{\ge0}}$ satisfying that
\[
\mathfrak{a}_m\cdot \mathfrak{a}_n\subset\mathfrak{a}_{n+m}
\]
for any $n,m\in\mathbb{Z}_{>0}$.
For any prime divisor $E$ over $X$, we set (cf.~\cite[Lemma 2.3]{JM})
\[
\mathrm{ord}_E(\mathfrak{a}_{\bullet}):=\lim_{m\to\infty}\frac{\mathrm{ord}_E(\mathfrak{a}_{m})}{m}=\inf_{m\ge0}\frac{\mathrm{ord}_E(\mathfrak{a}_{m})}{m}.
\]

By \cite[Corollary 2.16]{JM}, we can set 
\[
\mathrm{lct}(X,B;\mathfrak{a}_\bullet):=\lim_{m\to\infty}m\cdot\mathrm{lct}(X,B;\mathfrak{a}_m)=\inf_E\frac{A_{(X,B)}(E)}{\mathrm{ord}_E(\mathfrak{a}_{\bullet})}.
\]
\end{de}

$(X,B,L)$ is called a {\it polarized log pair} if $(X,B)$ is a log pair and $L$ is an ample $\mathbb{Q}$-line bundle.

\begin{de}
    Let $(X,B,L)$ be a polarized klt pair and take $r_0\in\mathbb{Z}_{>0}$ such that $r_0L$ is Cartier.
    For any $m\in\mathbb{Z}_{>0}$, we call $D$ an $mr_0${\it -basis type} divisor of $L$ if there exists a basis $\{D_i\}_{i=1}^{h^0(X,\mathcal{O}_X(mr_0L))}$ of $H^0(X,\mathcal{O}_X(mr_0L))$ such that $$D=\frac{1}{mr_0h^0(X,\mathcal{O}_X(mr_0L))}\sum_{i=1}^{h^0(X,\mathcal{O}_X(mr_0L))}D_i.$$
    $|L|_{mr_0\text{-basis}}$ denotes the set of all $mr_0$-basis type divisors of $L$.
    We set $$\delta_{r_0m}(X,B,L):=\inf_{D\in|L|_{mr_0\text{-basis}}}\mathrm{lct}(X,B;D)$$
    and $\delta(X,B,L):=\limsup_{m\to\infty}\delta_{r_0m}(X,B,L)$.
    We call $\delta(X,B,L)$ the $\delta${\it -invariant} and know by \cite[Theorem A]{BlJ} that $\delta(X,B,L)=\lim_{m\to\infty}\delta_{r_0m}(X,B,L)$.
\end{de}
Let $(X,B,L)$ be a polarized klt pair as above.
$|L|_{\mathbb{Q}}$ denotes the set of all effective $\mathbb{Q}$-divisors $\mathbb{Q}$-linearly equivalent to $L$.
Suppose that $r_0L$ is Cartier for some $r_0\in\mathbb{Z}_{>0}$.
We set for any prime divisor $E$ over $X$, 
\[
T_L(E)=\sup_{D\in |L|_{\mathbb{Q}}}\mathrm{ord}_E(D),\quad S_L(E)=\lim_{m\to \infty}\sup_{D\in |L|_{m\text{-basis}}}\mathrm{ord}_E(D).
\]
Indeed, we see that the above limit exists (cf.~\cite{BlJ}).
We set the $\alpha${\it -invariant} as $$\alpha(X,B,L):=\inf_E\frac{A_{(X,B)}(E)}{T_L(E)}.$$
On the other hand, we see that 
\[
\delta(X,B,L)=\inf_E\frac{A_{(X,B)}(E)}{S_L(E)}.
\]

\begin{de}[K-stability]
    Let $(X,B,L)$ be a polarized log pair of dimension $n$.
    A {\it normal semiample test configuration} $(\mathcal{X},\mathcal{L})$ for $(X,L)$ is defined as follows.
    \begin{enumerate}
    \item $\mathcal{X}$ is a normal variety with a $\mathbb{G}_m$-action and $\mathcal{L}$ is a semiample $\mathbb{G}_m$-linearized $\mathbb{Q}$-line bundle on $\mathcal{X}$. 
        \item There exists a proper surjective and $\mathbb{G}_m$-equivariant morphism $\pi\colon\mathcal{X}\to\mathbb{A}^1$, where $\mathbb{G}_m$ acts on $\mathbb{A}^1$ by multiplication.
        \item $(\pi^{-1}(1),\mathcal{L}|_{\pi^{-1}(1)})\cong (X,L)$.
    \end{enumerate}
   Let $(X_{\mathbb{A}^1},L_{\mathbb{A}^1})$ denote a semiample test configuration $(X\times\mathbb{A}^1,L\times\mathbb{A}^1)$ with the trivial $\mathbb{G}_m$-action on the first component $X$.
It is well-known that there exists another semiample test configuration $(\mathcal{Y},\mathcal{L}_{\mathcal{Y}})$ and two $\mathbb{G}_{m}$-equivariant birational morphism $\sigma\colon\mathcal{Y}\to\mathcal{X}$ and $\rho\colon\mathcal{Y}\to X_{\mathbb{A}^1}$ such that $\mathcal{L}_{\mathcal{Y}}=\sigma^*\mathcal{L}$ and their restrictions to $(\pi\circ\sigma)^{-1}(1)$ are nothing but the identity morphisms of $X$. 
Then, consider the $\mathbb{G}_m$-equivariant canonical compactification $(\overline{\mathcal{Y}},\overline{\mathcal{L}_{\mathcal{Y}}})\to\mathbb{P}^1$ of $(\mathcal{Y},\mathcal{L}_{\mathcal{Y}})\to\mathbb{A}^1$ such that the fiber over $\infty\in\mathbb{P}^1$ coincides with $(X,L)$ with the trivial $\mathbb{G}_m$-action.
Let $H$ be an arbitrary $\mathbb{R}$-line bundle over $X$ and let  $\mathcal{B}_{\mathcal{Y}}$ be the Zariski-closure of $B\times\mathbb{G}_m$ in $\overline{\mathcal{Y}}$.
Then we set the non-Archimedean Mabuchi functional and the non-Archimedean $\mathrm{J}^H$-functional as
\begin{align*}
M^{\mathrm{NA}}_B(\mathcal{X},\mathcal{L})&=(K_{\overline{\mathcal{Y}}/\mathbb{P}^1}+\mathcal{B}_{\mathcal{Y}}-\mathcal{Y}_0+\mathcal{Y}_{0,\mathrm{red}})\cdot\overline{\mathcal{L}_{\mathcal{Y}}}^n-\frac{n(K_X+B)\cdot L^{n-1}}{(n+1)L^n}\overline{\mathcal{L}_{\mathcal{Y}}}^{n+1} \\
    \mathcal{J}^{H,\mathrm
NA}(\mathcal{X},\mathcal{L})&=\overline{\rho^*H_{\mathbb{A}^1}}\cdot\overline{\mathcal{L}_{\mathcal{Y}}}^n-\frac{nH\cdot L^{n-1}}{(n+1)L^n}\overline{\mathcal{L}_{\mathcal{Y}}}^{n+1}.
\end{align*}
It is well-known that $M^{\mathrm{NA}}_B(\mathcal{X},\mathcal{L})$ and
$\mathcal{J}^{H,\mathrm
NA}(\mathcal{X},\mathcal{L})$ do not depend on the choice of $(\mathcal{Y},\mathcal{L}_{\mathcal{Y}})$ (cf.~\cite[\S7]{BHJ} and \cite{Hat2}).

We say that $(X,B,L)$ is {\it uniformly} K-{\it stable} (resp.~$(X,L)$ is {\it uniformly $\mathrm{J}^H$-stable}) if there exists $\epsilon>0$ such that  
$$M^{\mathrm{NA}}_B(\mathcal{X},\mathcal{L})\quad (\text{resp.~}\mathcal{J}^{H,\mathrm
NA}(\mathcal{X},\mathcal{L}))\ge\epsilon\mathcal{J}^{L,\mathrm
NA}(\mathcal{X},\mathcal{L})$$
for any normal semiample test configuration $(\mathcal{X},\mathcal{L})$ for $(X,L)$.
We say that $(X,B,L)$ is K{\it -semistable} (resp.~$(X,L)$ is $\mathrm{J}^H${\it -semistable}) if 
$$M^{\mathrm{NA}}_B(\mathcal{X},\mathcal{L})\quad (\text{resp.~}\mathcal{J}^{H,\mathrm
NA}(\mathcal{X},\mathcal{L}))\ge0$$
for any normal semiample test configuration $(\mathcal{X},\mathcal{L})$ for $(X,L)$.
It is well-known that $\mathcal{J}^{L,\mathrm
NA}(\mathcal{X},\mathcal{L})\ge0$ (cf.~\cite[Proposition 7.8]{BHJ}) and hence uniform K-stability implies K-semistability.
\end{de}

For J$^H$-semistability in the case when $H$ is ample, there exists a useful criterion.
We prepare the following notion.
\begin{de}
    Let $(X,L)$ be a polarized normal variety of dimension $n$ with an $\mathbb{R}$-line bundle $H$.
We say that $(X,L)$ is $\mathrm{J}^H${\it -nef} if for any subvariety $V\subset X$ of dimension $p$,    \[
\left(n\frac{H\cdot L^{n-1}}{L^n}L^p-pH\cdot L^{p-1}\right)\cdot V\ge0.
\]
We also say that $(X,L)$ is {\it uniformly} $\mathrm{J}^H${\it -positive} if there exists $\epsilon>0$ such that $(X,L)$ is $\mathrm{J}^{H-\epsilon L}$-nef.
\end{de}

\begin{thm}[{\cite[Theorem 3.2]{Hat}}]\label{thm-j-nef-j-st}
Let $(X,L)$ be a polarized normal variety with a nef $\mathbb{R}$-line bundle $H$.
Then $(X,L)$ is $\mathrm{J}^{H}$-nef if and only if $(X,L)$ is $\mathrm{J}^{H}$-semistable.

If $H$ is further ample, then $(X,L)$ is uniformly $\mathrm{J}^{H}$-positive if and only if $(X,L)$ is uniformly $\mathrm{J}^{H}$-stable.
\end{thm}

By using J-stability and the $\delta$-invariant, we can set the following key notion.

\begin{de}[Special K-stability, {\cite[Definition 3.10]{CM}}]
    Let $(X,B,L)$ be a polarized klt pair.
    If $(X,B,L)$ is uniformly J$^{\delta(X,B,L)L+K_X+B}$-stable and $\delta(X,B,L)L+K_X+B$ is ample, then we say that $(X,B,L)$ is {\it specially} K-{\it stable}.
    If $(X,B,L)$ is J$^{\delta(X,B,L)L+K_X+B}$-semistable and $\delta(X,B,L)L+K_X+B$ is nef, then we say that $(X,B,L)$ is {\it specially} K-{\it semistable}.
\end{de}

By \cite[Corollary 3.21]{CM}, specially K-stable log pairs are uniformly K-stable.

\subsection{Filtered linear series and good filtrations}

We collect some fundamental concepts of filtrations.

\begin{de}
    Let $R=\oplus_{m\in\mathbb{Z}_{\ge0}}R_m$ be a finitely generated graded $\mathbb{C}$-algebra such that $R_0=\mathbb{C}$.
    $F$ is called a ({\it decreasing, left-continuous and multiplicative}) {\it filtration} of $R$ if $F^{\lambda}R_m\subset R_m$ is a vector subspace for any $\lambda\in\mathbb{R}$ and $m\in\mathbb{Z}_{\ge0}$ and the following hold.
    \begin{itemize}
        \item $F^{\lambda}R_m\subset F^{\lambda'}R_m$ for any $\lambda>\lambda'$ and $F^{\lambda}R_m=\bigcap_{\lambda>\lambda'} F^{\lambda'}R_m$,
        \item $F^{\lambda}R_m\cdot F^{\lambda'}R_{m'}\subset F^{\lambda+\lambda'}R_{m+m'}$, and
        \item $F^0R_0=R_0$.
    \end{itemize}
    We say that this filtration is {\it linearly bounded} if there exist a positive real number $C$ and $m_0\in\mathbb{Z}_{>0}$ such that $F^{\lambda}R_m=0$ (resp.~$=R_m$) for any $m\ge m_0$ and $\lambda\ge Cm$ (resp.~$\le -Cm$).

Let $X$ be a normal projective variety and $L$ an ample line bundle on $X$.
    Let $F$ be a linearly bounded multiplicative filtration on $R=\oplus_{m\ge0}H^0(X,\mathcal{O}_X(mL))$.
    Then we define a graded subalgebra $FR^{(\lambda)}:=\oplus_{m\ge0}F^{m\lambda}H^0(X,\mathcal{O}_X(mL))\subset R$.
 We set $$\lambda_{\max}(F):=\limsup_{k\to \infty}\frac{\sup\{t\in\mathbb{R}|F^{t}R_k\ne0\}}{k}$$  
 and 
 \[
 \lambda_{\min}(F):=\inf\{\lambda\in\mathbb{R}|\mathrm{vol}(FR^{(\lambda)})>0\}.
 \]
 We set the {\it weight} $w_{F}(m)$ of $R_m$ with respect to $F$ as
    \[
    w_F(m)=\sum_{\lambda\in\mathbb{R}}\mathrm{dim}\left(F^{\lambda}R_m/\bigcup_{\lambda<\lambda'} F^{\lambda'}R_m\right)
    \]
    and we call $w_{F}$ the {\it weight function} of $F$.    
    Set 
    \begin{align*}
    S_{m}(F)&:=\frac{w_F(m)}{mh^0(X,\mathcal{O}_X(mL))},\quad\text{and}\\
    S(F)&:=\lambda_{\min}(F)+\frac{1}{(L^n)}\int^{\lambda_{\max}(F)}_{\lambda_{\min}(F)}\mathrm{vol}(FR^{(\lambda)})d\lambda.
\end{align*}
It is well-known that $S(F)=\lim_{m\to\infty}S_m(F)$ (cf.~\cite[Corollary 2.12]{BlJ}). 
\end{de}
There are following two fundamental examples.
\begin{ex}\label{ex--restricted--filtration}
    Let $X$ be a proper variety with an ample line bundle $L$ on $X$.
    Let $R:=\oplus H^0(X,\mathcal{O}_X(mL))$ and $D$ a closed subvariety of $X$.
    Let $F$ be a linearly bounded multiplicative filtration on $R$.
    Then we set a filtration $F|_D$ on $\oplus H^0(D,\mathcal{O}_D(mL|_D))$ as
    \[
    F|_D^\lambda H^0(D,\mathcal{O}_D(mL|_D)):=\mathrm{Image}(F^\lambda H^0(X,\mathcal{O}_X(mL))\to H^0(D,\mathcal{O}_D(mL|_D)))
    \]
    for any $\lambda\in\mathbb{R}$ and $m\in\mathbb{Z}_{\ge0}$.
    Then we call $F|_D$ the {\it restricted filtration} on $D$.
    We can check that $F|_D$ is multiplicative and linearly bounded.
\end{ex}

\begin{ex}
    Let $F$ be a linearly bounded multiplicative filtration of $R=\oplus_{m\ge0}R_m$.
    If we set for any $\lambda\in\mathbb{R}$,  $$F_{\mathbb{Z}}^{\lambda}R_m:=F^{\lceil\lambda\rceil}R_m,$$
    then we see that $F_{\mathbb{Z}}$ is also a linearly bounded multiplicative filtration.
\end{ex}

\begin{de}[Good filtrations]
    Let $X$ be a normal projective variety of dimension $n$ and $L$ an ample line bundle on $X$.
    Let $F$ be a linearly bounded multiplicative filtration on $R=\oplus_{m\ge0}H^0(X,\mathcal{O}_X(mL))$.
    If there exist $a_0,a_1\in\mathbb{R}$ and $C\in\mathbb{R}_{>0}$ such that 
    \[
    |w_F(m)-a_0m^{n+1}-a_1m^n|<Cm^{n-2}
    \]
    for any $m\in\mathbb{Z}_{>0}$, then we say that $F$ is a {\it good filtration}.
\end{de}

We note that if $F=F_{\mathbb{Z}}$ and $\oplus F^\lambda R_m$ is a finitely generated bigraded algebra, then $F$ is good. 

The following notion is first introduced by \cite{XZ}.

\begin{de}
    Let $(X,B,L)$ be a polarized klt pair. 
Take $r\in\mathbb{Z}_{>0}$ such that $rL$ is Cartier and a linearly bounded filtration $F$ on $R:=\oplus_{m\ge0}R_m$, where $R_m=H^0(X,\mathcal{O}_X(mrL))$.
Now, we set the base ideal of $F^{\lambda}R_m$ as
\[
I_{m,\lambda}(F):=\mathrm{Image}(F^{\lambda}R_m\otimes\mathcal{O}_{X}(-mrL)\to\mathcal{O}_X)
\]
for $\lambda\in\mathbb{R}$.
Set $I_\bullet^{(\lambda)}(F):=\{I_{m,m\lambda}(F)\}_{m\ge0}$ as a graded sequence of ideals.
For any $\delta>0$, we set the $\delta$-{\it lc slope} of $F$ as
\[
\mu_\delta(F):=\sup\left\{\lambda\in\mathbb{R}\,\bigg|\,\mathrm{lct}(X,B;I_{\bullet}^{(\lambda)}(F))\ge\frac{\delta}{r}\right\}.
\]
Furthermore, we set
\[
\beta_{\delta}(F):=\frac{\mu_\delta(F)-S(F)}{r}.
\]
\end{de}

We remark that the following holds.
\begin{thm}[{\cite[Proposition 4.5]{XZ}}]\label{thm--xz--4.5}
    Let $(X,B,L)$ be a polarized klt pair. 
Take $r\in\mathbb{Z}_{>0}$ such that $rL$ is Cartier and set $R:=R(X,rL)$.
Then,
\[
\delta(X,\Delta,L)=\sup\{\delta>0\,|\,\beta_{\delta}(F)\ge 0\text{ for any linearly bounded filtration }F\text{ on $R$}\}.
\]
\end{thm}

We remark that Xu and Zhuang showed Theorem \ref{thm--xz--4.5} only for log Fano pairs but their proof also works for Theorem \ref{thm--xz--4.5} in the same way.
We also remark the following useful lemma, which also holds for general polarized klt pairs.

\begin{lem}[{\cite[Lemma 4.13]{XZ}}]\label{lem--xz--4.13}
Let $(X,B,L)$ be a polarized klt pair. 
Take $r\in\mathbb{Z}_{>0}$ such that $rL$ is Cartier and a linearly bounded filtration $F$ on $R:=\oplus_{m\ge0}R_m$, where $R_m=H^0(X,\mathcal{O}_X(mrL))$.
For any real numbers  $s,\epsilon\in(0,1)$, it holds that
\[
\mu_{1+(1-\epsilon)s}(F)\ge s\cdot\mu_{\epsilon^{-1}}(F)+(1-s)\mu_{1}(F).
\]
\end{lem}

\begin{de}[Donaldson-Futaki invariant and J-functional for filtrations]\label{de--j--filt}
    Let $(X,L)$ be a polarized variety of dimension $n$ with $L$ a line bundle and let $F$ be a multiplicative filtration on $R$, where $R=\oplus_{m\ge0}R_m$ and $R_m=H^0(X,\mathcal{O}_X(mL))$.
    Let $w_F(m)$ be the weight function of $F$.
Then we see that (cf.~\cite[(2)]{CM})
\[
\lim_{m\to\infty}\frac{w_F(m)}{m^{n+1}}
\]
exists and write this by $b_0$.
On the other hand, it is well-known that $\chi(X,\mathcal{O}_X(mL))$ is a polynomial of degree $n$.
We denote this by $a_0m^n+a_1m^{n-1}+O(m^{n-2})$.

Take $H$ an ample $\mathbb{Q}$-line bundle on $X$.
Take a sufficiently divisible $r\in\mathbb{Z}_{>0}$ such that $rH$ is very ample.
Let $D\in |rH|$ be a very general member and let $F|_D$ be the restriction (cf.~Example \ref{ex--restricted--filtration}).
Let $\tilde{a}_0:=\frac{H\cdot L^{n-1}}{(n-1)!}$ and $\tilde{b}_0:=\lim_{m\to\infty}\frac{w_{F|_D}(m)}{rm^n}$.
By \cite[Lemma 2.20]{CM}, we see that the following value
\[
\frac{\tilde{b}_0a_0-\tilde{a}_0b_0}{a_0^2}
\]
is independent from the choice of $r$ and very general $D$.
We denote this by $\mathcal{J}^{H,\mathrm{NA}}(F)$ and call the $\mathrm{J}^H$-functional of $F$.
For an arbitrary $\mathbb{Q}$-line bundle $T$ on $X$, there exist two ample line bundles $H_1$ and $H_2$ such that $T=H_1-H_2$.
It is easy to see that $\mathcal{J}^{H,\mathrm{NA}}(F)$ is linear with respect to $H$.
We set $\mathcal{J}^{T,\mathrm{NA}}(F):=\mathcal{J}^{H_1,\mathrm{NA}}(F)-\mathcal{J}^{H_2,\mathrm{NA}}(F)$.
$\mathcal{J}^{T,\mathrm{NA}}(F)$ is independent from the choice of $H_1$ and $H_2$.

If $F$ is a good filtration and $w_F(m)=b_0m^{n+1}+b_1m^n+O(m^{n-1})$, then we set the Donaldson-Futaki invariant of $F$ as
\[
\mathrm{DF}(F)=2\frac{b_0a_1-a_0b_1}{a_0^2}.
\]
\end{de}

\begin{prop}[{\cite[Lemma 2.20]{CM}}]\label{prop--filt--j}
Let $(X,L)$ be a polarized variety with $L$ a line bundle and let $H$ be a nef $\mathbb{Q}$-line bundle.
If $(X,L)$ is $\mathrm{J}^H$-semistable, then 
\[
\mathcal{J}^{H,\mathrm{NA}}(F)\ge0
\]
for any multiplicative linearly bounded filtration $F$ on $\oplus_{m\ge0} H^0(X,\mathcal{O}_X(mL))$.
\end{prop}

\begin{proof}
    Let $F$ be a multiplicative linearly bounded filtration on $\oplus_{m\ge0} H^0(X,\mathcal{O}_X(mL))$.
   By the fact that $\mathcal{J}^{H,\text{NA}}(F)$ is linear with respect to $H$, we may assume that $H$ is very ample and take a very general member $D\in|H|$.
    By \cite[Lemma 2.20]{CM}, we see that
    \[
\mathcal{J}^{H,\mathrm{NA}}(F_{\mathbb{Z}})\ge0.
\]
By \cite[Corollary 2.12]{BlJ}, 
$S(F)=S(F_{\mathbb{Z}})$.
On the other hand, $F_{\mathbb{Z}}|_D=(F|_D)_\mathbb{Z}$.
Since 
\[
\mathcal{J}^{H,\mathrm{NA}}(F)=n\left(S(F|_D)-S(F)\right),
\]
we have $\mathcal{J}^{H,\mathrm{NA}}(F_{\mathbb{Z}})=\mathcal{J}^{H,\mathrm{NA}}(F)$, which completes the proof.
\end{proof}

\subsection{Polarized log family}
In this subsection, we discuss the following concept.
\begin{de}\label{de--polarized--family}
    Let $\pi\colon X\to S$ be a proper flat morphism of normal varieties such that $\pi_*\mathcal{O}_X\cong\mathcal{O}_S$, $\Delta$ an effective $\mathbb{Q}$-Weil divisor on $X$ and $L$ a $\pi$-ample $\mathbb{Q}$-line bundle on $X$.
    We say that $\pi\colon(X,\Delta,L)\to S$ is a {\it polarized log family} if any fiber $X_s$ over $s\in S$ is normal and no irreducible component of $\Delta$ contains some fiber $X_s$.
    If $K_{X/S}+\Delta$ is $\mathbb{Q}$-Cartier, then we say that $\pi\colon(X,\Delta,L)\to S$ is {\it $\mathbb{Q}$-Gorenstein}.
We say that $\pi$ is of relative dimension $n$ if $\mathrm{dim}\,X_s=n$ for general $s\in S$. 
\end{de}

A polarized log $\mathbb{Q}$-Gorenstein family $\pi\colon(X,\Delta,L)\to S$ satisfies the condition of (\ref{Notations--(vi)}) in Notations and convention.
Hence, for any morphism $g\colon T\to S$ from a normal variety, we can set $\Delta_T$ as (\ref{Notations--(vi)}).
To state Theorem \ref{thm--main--2}, we prepare the following notion.
\begin{de}[Maximal variation]\label{de--maximal-variation}
   Let $\pi\colon(X,\Delta,L)\to S$ be a polarized log family.
   Suppose that for any irreducible curve $C\subset S$ containing a general point of $S$ and two general distinct closed points $p,q\in C$, $(X_p,\Delta_p)$ and $(X_q,\Delta_q)$ are not isomorphic.
   Then we say that $\pi$ has {\it maximal variation}.
\end{de}
Next, we show that if some fiber of a polarized log $\mathbb{Q}$-Gorenstein family is specially K-stable, then so are very general fibers.
We make use of this assertion to show Theorem \ref{thm--main--2}.
To prove this, we first show the following on J-semistability.

\begin{prop}\label{prop--generic-j-ss}
    Let $\pi\colon(X,\Delta,L)\to S$ be a polarized log family and let $H$ be a $\mathbb{Q}$-line bundle on $X$. Suppose that there exists a closed point $s_0\in S$ such that $H_{s_0}$ is nef and $(X_{s_0},L_{s_0})$ is $\mathrm{J}^{H_{s_0}}$-semistable.
    Then for any very general point $s\in S$, $(X_{s},L_{s})$ is $\mathrm{J}^{H_{s}}$-semistable.
\end{prop}

\begin{proof}
    It is well-known that if $H_{s_0}$ is nef, then so is $H_{s}$ for any very general point $s\in S$.
On the other hand, we deal with J$^{H_s}$-nefness. 
Consider the Hilbert scheme $\mathrm{Hilb}_{X/S}$.
Recall that for any $S$-scheme $T$, the set of $T$-valued points of $\mathrm{Hilb}_{X/S}$ is the set of closed subschemes of $X\times_{S}T$ flat over $T$.
It is well-known that $\mathrm{Hilb}_{X/S}$ has countably many connected components and each of them is proper over $S$.
Here, we claim that for very general point $s\in S$ and any $p$-dimensional closed subvariety $V\subset X_s$, 
\[
\left(n\frac{H_s\cdot L_s^{n-1}}{L_s^n}L_s-pH_s\right)\cdot L_s^{p-1}\cdot V\ge0.
\]
Indeed, for any closed subvariety $V\subset X_s$ of dimension $p$, there exists a connected component $\mathcal{H}\subset \mathrm{Hilb}_{X/S}$ containing the point corresponding to $V$.
Since $s$ is very general, $\mathcal{H}\to S$ is surjective.
 Then \[
\left(n\frac{H_s\cdot L_s^{n-1}}{L_s^n}L_s-pH_s\right)\cdot L_s^{p-1}\cdot V=\left(n\frac{H_{s_0}\cdot L_{s_0}^{n-1}}{L_{s_0}^n}L_{s_0}-pH_{s_0}\right)\cdot L_{s_0}^{p-1}\cdot V'\ge0
\]
by the J$^{H_{s_0}}$-nefness of $(X_{s_0},L_{s_0})$, where $V'\subset X_{s_0}$ is a $p$-dimensional closed subscheme whose corresponding point in $\mathrm{Hilb}_{X/S}$ is contained in $\mathcal{H}$.
This means that for any very general closed point $s\in S$, $(X_{s},L_{s})$ is $\mathrm{J}^{H_{s}}$-nef and $H_s$ is nef. By Theorem \ref{thm-j-nef-j-st}, we obtain the assertion.
\end{proof}

\begin{cor}\label{cor--verygeneral-special-kst}
    Let $\pi\colon(X,\Delta,L)\to S$ be a polarized log $\mathbb{Q}$-Gorenstein family. 
    Suppose that there exists a closed point $s_0\in S$ such that $(X_{s_0},\Delta_{s_0},L_{s_0})$ is specially {\rm K}-stable.
    Then, 
    there exist positive rational numbers $\lambda$ and $\epsilon$ such that for any very general point $s\in S$, $\delta{(X_s,\Delta_s,L_s)}\ge\lambda +\epsilon$ and $(X_s,L_s)$ is $\mathrm{J}^{K_{X_s}+\Delta_s+\lambda L_s}$-semistable.
    In particular, $(X_{s},\Delta_s,L_{s})$ is specially {\rm K}-stable for any very general closed point $s\in S$.
\end{cor}

\begin{proof}
We know that the correspondence
\[
S\ni s\mapsto \delta(X_{\bar{s}},\Delta_{\bar{s}},L_{\bar{s}})
\] 
is lower-semicontinuous by \cite[Theorem 6.6]{BL}.
Thus, for any sufficiently small $\epsilon\in\mathbb{Q}_{>0}$, there exists a non-empty open subset $U\subset S$ such that $\delta(X_{\bar{s}},\Delta_{\bar{s}},L_{\bar{s}})\ge \delta(X_{s_0},\Delta_{s_0},L_{s_0})-\epsilon$ for any geometric point $\bar{s}\in U$.
Choose $\epsilon$ small enough such that there exists $\lambda\in\mathbb{Q}_{>0}$ such that $(X_{s_0},L_{s_0})$ is uniformly $\mathrm{J}^{K_{X_{s_0}}+\Delta_{s_0}+\lambda L_{s_0}}$-stable and $\lambda+2\epsilon\le \delta(X_{s_0},\Delta_{s_0},L_{s_0})$.
Then we see that $\delta{(X_{\bar{s}},\Delta_{\bar{s}},L_{\bar{s}})}\ge\lambda +\epsilon$ for any geometric point $\bar{s}\in U$.
On the other hand, $(X_s,L_s)$ is $\mathrm{J}^{K_{X_s}+\Delta_s+\lambda L_s}$-semistable by Proposition \ref{prop--generic-j-ss} for any very general closed point $s\in S$.
\end{proof}

\section{CM line bundle}

In this section, we discuss the CM line bundle. 
\subsection{CM line bundle and J-line bundle}
First, we explain how to define the log CM line bundle for a polarized log $\mathbb{Q}$-Gorenstein family.

\begin{de}
    Let $S$ be a Noetherian scheme and let $E$ be a vector bundle over $S$.
Let $X\subset \mathbb{P}_S(E)$ be a closed subscheme such that $\mathcal{O}_X$ is a perfect complex as an $\mathcal{O}_{\mathbb{P}_S(E)}$-module.
Suppose that the generic fiber of the canonical morphism $\pi\colon X\to S$ is of dimension $n$.
We say that $\mathcal{O}_X$ satisfies the condition $Q_{(r)}$ if the following holds:
\begin{enumerate}
    \item for each point $y\in Y$ of depth $0$,
    \[
    \mathrm{dim}((\mathrm{Supp}(\mathcal{O}_{X}))_y)\le r,
    \]
     \item for each point $y\in Y$ of depth $1$,
    \[
    \mathrm{dim}((\mathrm{Supp}(\mathcal{O}_{X}))_y)\le r+1.
    \]
\end{enumerate}
Let $L:=\mathcal{O}_{\mathbb{P}_S(E)}(1)|_X$.
Assume now that $\mathcal{O}_X$ satisfies the condition $Q_{(r)}$. 
For any sufficiently large $m\in\mathbb{Z}_{>0}$,
    consider the {\it Knudsen-Mumford expansion} \cite[Theorem 4]{KnMu} (cf.~\cite[Lemma 5.8]{GIT})
    \[
    \mathrm{det}(\pi_*\mathcal{O}_X(mL))\cong\bigotimes_{i=0}^{n+1}\mathcal{M}_i^{\otimes\binom{m}{i}},
    \]
    where $\mathcal{M}_i$ is a uniquely determined line bundle on $S$ for $i=0,\ldots,n+1$.
     It is well-known that the Knudsen-Mumford expansion is compatible with base changes.
    More precisely, for any morphism $g\colon T\to S$, consider the Knudsen-Mumford expansion 
    \[
    \mathrm{det}(\pi_{T*}\mathcal{O}_{X_T}(mL_T))\cong\bigotimes_{i=0}^{n+1}\mathcal{N}_i^{\otimes\binom{m}{i}},
    \]
    where $\pi_T$, $X_T$ and $L_T$ are the base changes by $T\to S$.
    Then $\mathcal{N}_i=g^*\mathcal{M}_i$ for $0\le i\le n+1$ (cf.~\cite[Lemma 3.5]{CP}).
\end{de}

\begin{de}[CM-line bundle]\label{de-cm-line-bundle}
    Let $\pi\colon(X,\Delta,L)\to S$ be a polarized log family of relative dimension $n$.
    Take $r\in\mathbb{Z}_{>0}$ such that $rL$ is a $\pi$-very ample line bundle.
    Then, $X\subset\mathbb{P}_S(\pi_*\mathcal{O}_X(mrL))$.
    We set the {\it log CM line bundle} of $\pi$ as
    \[
    \lambda_{\mathrm{CM},\pi}:=\pi_*(\mu_LL^{n+1}+(n+1)L^n\cdot(K_{X/S}+\Delta)),
    \]
    where $\mu_L:=n\frac{-(K_{X_t}+\Delta_t)\cdot L_t^{n-1}}{L_t^n}$ for general closed point $t\in S$ and $\pi_*(L^n\cdot D)$ is defined to be a $\mathbb{Q}$-divisor unique up to $\mathbb{Q}$-linear equivalence on $S$ for any $\mathbb{Q}$-Cartier $\mathbb{Q}$-divisor $D$ on $X$ as follows.
    Suppose that $mL$ is relatively very ample over $S$.
    Take a line bundle $M$ on $S$ such that $N:=mL+\pi^*M$ is very ample.
    By the Bertini theorem and \cite[II, Exercise 8.2]{Ha}, we see that there exist a positive integer $l$ and $D_1,\ldots,D_{n}\in|lN|$ such that $Y:=D_1\cap D_2\cap\ldots\cap D_n$ is normal, irreducible and finite over $S$.
    We may further assume that $Y\not\subset \mathrm{Supp}\,D$ and then we can define $D\cap Y$ as a $\mathbb{Q}$-Cartier $\mathbb{Q}$-divisor on $Y$.
    Then we set a $\mathbb{Q}$-Weil divisor
    \[
    \pi_*(L^n\cdot D):=\frac{1}{(ml)^n}\pi_*(D\cap Y)-nm(D_t\cdot L_{t}^{n-1})M.
    \]
    Then we see the following.

    \begin{prop}\label{prop--base-change-cm}
        Let $\pi\colon(X,\Delta,L)\to S$
be a log $\mathbb{Q}$-Gorenstein polarized family with a $\mathbb{Q}$-Cartier $\mathbb{Q}$-divisor $D$ on $X$ as above.
Then 
$\pi_*(L^n\cdot D)$ is a $\mathbb{Q}$-Cartier $\mathbb{Q}$-divisor on $S$ uniquely determined up to $\mathbb{Q}$-linear equivalence independent from the choices of $m$, $l$ and $D_1,\ldots,D_n$.
In particular, $\lambda_{\mathrm{CM},\pi}$ is a well-defined $\mathbb{Q}$-line bundle on $S$ and for any morphism $g\colon T\to S$ from a normal variety, it holds that $\lambda_{\mathrm{CM},\pi_T}=g^*(\lambda_{\mathrm{CM},\pi})$.
    \end{prop}
\begin{proof}
We note that for any $\mathbb{Q}$-Cartier $\mathbb{Q}$-divisors $E_1$ and $E_2$ on $X$, we have that $$\pi_*(Y\cap (E_1+E_2))\sim_{\mathbb{Q}}\pi_*(Y\cap E_1)+\pi_*(Y\cap E_2)$$ and if $E_1\sim_{\mathbb{Q}}E_2$, then $\pi_*(Y\cap E_1)\sim_{\mathbb{Q}}\pi_*(Y\cap E_2)$.
Thus, we may replace $D$ by a very ample line bundle on $X$.
We may further assume that $D$ is an effective normal Cartier divisor on $X$ such that every fiber of $\pi|_D\colon D\to S$ is equidimensional and of dimension $n-1$ by the Bertini theorem and \cite[II, Exercise 8.2]{Ha}.
Then we see that $\mathcal{O}_D$ satisfies the assumption of \cite[Theorem 4]{KnMu}.
Indeed, it suffices to show that $\mathcal{O}_D$ is a perfect complex of $\mathbb{P}_S(E)$, where $E=\pi_*(\mathcal{O}_X(rL))$, since the condition $Q_{(r)}$ is satisfied (see \cite[p.~50]{KnMu} and \cite[Lemma A.1]{CP}).
For this, it is enough to show that $\mathcal{O}_D$ is an $\mathcal{O}_{\mathbb{P}_S(E)}$-module of finite Tor-dimension (cf.~\cite[p.~111]{GIT}).
By the fact that $D$ is a Cartier divisor of $X$ and \cite[Lemma 5.8]{GIT}, we know that $\mathcal{O}_D$ is also of finite Tor-dimension.
Thus, we may apply \cite[Theorem 4]{KnMu} for $\pi|_D$ and there exists the Knudsen-Mumford expansion \[
    \mathrm{det}((\pi|_D)_*\mathcal{O}_D(mrL|_D))\cong\bigotimes_{i=0}^{n}\mathcal{N}_i^{\otimes\binom{m}{i}}.
    \]
   We assert that $\pi_*(D\cdot L^n)=(\pi|_D)_*(L|_D^n)=\frac{1}{r^n}\mathcal{N}_n$.
   Indeed, take a big open subset $S^\circ\subset S$ such that $S^\circ$ is smooth and $D$ is flat over $S^\circ$.
   Over $S^\circ$,
   we have that 
   $$\pi_*(D\cdot L^n)|_{S^\circ}=(\pi|_D)_*(L|_D^n)|_{S^\circ}=\frac{1}{r^n}\mathcal{N}_n|_{S^\circ}$$
   by \cite[Lemmas A.1, A.2]{CP}.
    Hence, $\pi_*(D\cdot L^n)\sim_{\mathbb{Q}}\frac{1}{r^n}\mathcal{N}_n$ is $\mathbb{Q}$-Cartier and for any morphism $g\colon T\to S$ from a normal variety, we see that $g^*(\pi_*(L^{n}\cdot D))=\pi_{T*}(L_T^{n}\cdot D_T)$.
    It follows from this that $\lambda_{\mathrm{CM},\pi_T}=g^*(\lambda_{\mathrm{CM},\pi})$.
    We complete the proof.
    \end{proof}

\end{de}

\begin{de}[J-line bundle]\label{de--j--line}
Let $\pi\colon(X,\Delta,L)\to S$ be a polarized log pair with an $\mathbb{R}$-line bundle $H$ on $X$.
We set the $\mathrm{J}^H${\it -line bundle} with respect to $H$ as
\[
\lambda_{\mathrm{J},\pi,H}:=\pi_*\left((n+1)L^n\cdot H-n\frac{H_t\cdot L_t^{n-1}}{L_t^n}L^{n+1}\right),
\]
where $t\in S$ is a general closed point.
As Proposition \ref{prop--base-change-cm}, we have that for any morphism $g\colon T\to S$ from a normal variety, $g^*(\lambda_{\mathrm{J},\pi,H})=\lambda_{\mathrm{J},\pi_T,H_T}$.
\end{de}

\begin{de}[CM-degree and J-degree]
Let $\pi\colon(X,\Delta,L)\to C$ be a polarized log family of relative dimension $n$ with $C$ a proper smooth curve.    
Let $t$ be a closed point of $C$ and $v:=(L_t^n)$.
Then we set the CM degree as 
\[
\mathrm{CM}((X,\Delta,L)/C):=\frac{1}{(n+1)v}\mathrm{deg}_C\lambda_{\mathrm{CM},\pi}|_C.
\]

On the other hand, let $H$ be an $\mathbb{R}$-line bundle on $X$.
We set the $\mathrm{J}^H$-degree as
\[
\mathcal{J}^H((X,L)/C):=\frac{1}{(n+1)v}\mathrm{deg}_C\lambda_{\mathrm{J},\pi,H}|_C.
\]
\end{de}
We note that
\begin{align*}
\mathrm{CM}((X,\Delta,L)/C)&=v^{-1}\left((K_{X/S}+\Delta)\cdot L^{n}-\frac{n(K_{X_t}+\Delta_t)\cdot L_t^{n-1}}{(n+1)L_t^n}L^{n+1}\right),\\
\mathcal{J}^H((X,L)/C)&=v^{-1}\left(H\cdot L^{n}-\frac{nH_t\cdot L_t^{n-1}}{(n+1)L_t^n}L^{n+1}\right).
\end{align*}

Next, we explain a relationship between the CM degree and the Harder--Narasimhan filtration.

\begin{de}[Harder-Narasimhan filtration]\label{de-hn}
Let $C$ be a proper smooth curve.
For any locally free sheaf $E$ on $C$, we set the {\it slope} of $E$ as 
    \[
    \mu(E):=\frac{\mathrm{deg}_C\,E}{\mathrm{rank}\,E}.
    \]
    We say that $E$ is {\it semistable} if $\mu(E)\ge\mu(F)$ 
    for any nonzero subsheaf $ F\subsetneq E$.
It is well-known that there exists the unique sequence (cf.~\cite[Theorem 1.3.4]{HL})
\[
0=E_0\subsetneq E_1\subsetneq E_2\subsetneq\ldots \subsetneq E_{k-1}\subsetneq E_k=E
\]
such that $E_i/E_{i-1}$ is a semistable locally free sheaf and $\mu_i:=\mu(E_i/E_{i-1})$ satisfies that $\mu_i>\mu_{i+1}$ for $1\le i\le k$.
We denote $\mu_{\min}:=\mu_{k}$ and call this the {\it minimal slope} of $E$.
For any $\lambda\in\mathbb{R}$, we set 
$\mathscr{F}_{\mathrm{HN}}^{\lambda}E$ as the union of subsheaves of minimal slope at least $\lambda$.
We call $\mathscr{F}_{\mathrm{HN}}$ the {\it Harder-Narasimhan filtration} of $E$.
\end{de}

\begin{de}
    Let $\pi\colon(X,\Delta,L)\to C$ be a polarized log family of relative dimension $n$, where $C$ is a proper smooth curve.
    Take $r\in\mathbb{Z}_{>0}$ such that $rL$ is Cartier.
    Take $s\in C$ such that $X_s$ is normal.
    Let $\mathcal{R}_m:=\pi_*\mathcal{O}_X(mrL)$ and $R_m:=H^0(X_s,\mathcal{O}_{X_s}(mrL_s))$ for any $m\in\mathbb{Z}_{\ge0}$.
    We set the Harder-Narasimhan filtration $\mathscr{F}_{\mathrm{HN}}$ on $\mathcal{R}_m$ as Definition \ref{de-hn}.
   We define a filtration $F_{\mathrm{HN}}$ as
\[
F^\lambda_{\mathrm{HN}}R_m:=\mathrm{Image}(\mathscr{F}_{\mathrm{HN}}^{\lambda}\mathcal{R}_m\subset \mathcal{R}_m\to R_m).
\]
It is known by \cite[Lemma-Definition 2.26]{XZ} that the filtration $\mathscr{F}_{\mathrm{HN}}$ is linearly bounded and multiplicative and so is $F_{\mathrm{HN}}$. 
We call $F_{\mathrm{HN}}$ the {\it induced} filtration of $R=\oplus_{m\ge0} R_m$ (see \cite[\S2.8]{XZ}). 
\end{de}

The following is important to calculate the CM degree or the J-degree.

\begin{prop}\label{prop-cm-df-good}
    Let $\pi\colon(X,\Delta,L)\to C$ be a polarized log pair of relative dimension $n$, where $C$ is a proper smooth curve. 
     Let $s\in C$ be a closed point.
    Take $r\in\mathbb{Z}_{>0}$ such that $rL$ is Cartier and define the induced filtration $F_{\mathrm{HN}}$ on $R=\oplus_{m\ge0}H^0(X_s,\mathcal{O}_{X_s}(mrL_s))$.
    Then $F_{\mathrm{HN}}$ is a good filtration and $$S(F_{\mathrm{HN}})=r\frac{L^{n+1}}{(n+1)L_s^n}.$$
\end{prop}

\begin{proof}
    By definition (cf.~\cite[Prop.~4.6]{XZ}), we see that 
    \[
    S_m(F_{\mathrm{HN}})=\frac{\mathrm{deg}_C\,\pi_*\mathcal{O}_X(mrL)}{mh^0(X_s,\mathcal{O}_{X_s}(mrL_s))}.
    \]
    Furthermore, let $g(C)$ be the genus of $C$.
    By the Riemann-Roch theorem, 
    \begin{equation}
        \mathrm{deg}_C\,\pi_*\mathcal{O}_X(mrL)=h^0(X_s,\mathcal{O}_{X_s}(mrL_s))(g(C)-1)+\chi(X,\mathcal{O}_X(mrL)). \label{eq-riemann-roch}
    \end{equation}
    By this, we see that $F_{\mathrm{HN}}$ is a good filtration.
Note that 
\begin{align*}
    h^0(X_s,\mathcal{O}_{X_s}(mrL_s))&=\frac{(mr)^n}{n!}L_s^{n}+O(m^{n-1}),\\
    \chi(X,\mathcal{O}_X(mrL))&=\frac{(mr)^{n+1}}{(n+1)!}L^{n+1}+O(m^{n}).
\end{align*}
Thus we have the second assertion by $\lim_{m\to\infty}S_m(F_{\mathrm{HN}})=S(F_{\mathrm{HN}})$.
\end{proof}

As \cite[Corollary 3.9]{CM}, we obtain that the CM degree is nonnegative when a fiber is smooth and admits a unique cscK metric. 

\begin{cor}
    Notations as above.
    Suppose further that $\Delta=0$ and $X_s$ is a smooth variety with a cscK metric in $\mathrm{c}_1(L_s)$ and $\mathrm{Aut}(X_s,L_s)$ discrete for some $s\in C$.
    Then $\mathrm{CM}((X,0,L)/C)\ge0$.
\end{cor}

\begin{proof}
    Let $\chi(X,\mathcal{O}_X(mrL))=b_0m^{n+1}+b_1m^n+O(m^{n-1})$ and $h^0(X_s,\mathcal{O}_{X_s}(mrL_s))=a_0m^n+a_1m^{n-1}+O(m^{n-2})$.
    We see by (\ref{eq-riemann-roch}) that
    \[
    \mathrm{deg}_C\,\pi_*\mathcal{O}_X(mrL)=b_0m^{n+1}+(b_1+a_0(g(C)-1))m^n+O(m^{n-1}).
    \]
    Thus, we see that $\mathrm{CM}((X,0,L)/C)=\mathrm{DF}(F_{\mathrm{HN}})$.
    On the other hand, we see that $(X_s,L_s)$ is asymptotically Chow stable by \cite{Dn}.
    Thus, we see as \cite[Theorem 2.18]{CM} that
    \[
    \mathrm{DF}(F_{\mathrm{HN}})\ge0.
    \]
    We complete the proof.
\end{proof}

\subsection{Moduli of uniformly adiabatically K-stable klt-trivial fibrations over curves and the CM line bundle}

Let $f\colon X\to C$ be a morphism of normal varieties such that $f_*\mathcal{O}_X\cong\mathcal{O}_C$ and suppose that $C$ is a proper smooth curve.
We say that $f\colon (X,\Delta)\to C$ is a {\it klt--trivial fibration} over a curve if $K_{X}+\Delta\sim_{\mathbb{Q},C}0$ and $(X,\Delta)$ is klt.
Then we set the {\it discriminant divisor} $B:=\sum_{P\in C}(1-\mathrm{lct}(X,\Delta;f^*P))P$.
Let $M$ be a $\mathbb{Q}$-divisor on $C$ such that
\[
K_X+\Delta\sim_{\mathbb{Q}}f^*(K_C+M+B).
\]
We call $M$ the {\it moduli divisor}.
Take an $f$-ample $\mathbb{Q}$-Cartier Weil divisor $A$.
We say that $f\colon (X,\Delta,A)\to C$ is {\it uniformly adiabatically} K{\it -stable} if 
there exist  positive constants $\epsilon_0$ and $\delta$ such that 
\[
M^{\mathrm{NA}}_\Delta(\mathcal{X},\mathcal{M})\ge\delta\mathcal{J}^{\epsilon A+L,\mathrm{NA}}(\mathcal{X},\mathcal{M})
\]
for any $\epsilon\in(0,\epsilon_0)$ and normal semiample test configuration $(\mathcal{X},\mathcal{M})$ for $(X,\epsilon A+L)$, where $L$ is a fiber of $f$.
It is known by \cite[Theorem 1.1]{Hat} that $f$ is uniformly adiabatically K-stable if and only if one of the following holds.
\begin{itemize}
    \item
    $-K_X-\Delta$ is nef but not numerically trivial and $ \delta(C,B,-K_C-M-B)>1$, or
    \item $K_X+\Delta$ is nef.
\end{itemize}

Let $\mathfrak{Z}_{d,v,u}$ be the following subset for $d\in\mathbb{Z}_{>0}$, $u\in\mathbb{Q}_{>0}$ and $v\in\mathbb{Q}_{>0}$:
\[
\left\{
f\colon(X,0,A)\to \mathbb{P}^1
\;\middle|
\begin{array}{rl}
(i)&\text{$f$ is a uniformly adiabatically K-stable klt--trivial}\\
&\text{fibration over $\mathbb{P}^1$ with $\mathrm{dim}\,X=d$,}\\
(ii)&\text{$A$ is an $f$-ample line bundle such that}\\
&\text{$K_{X}\cdot A^{d-1}=-uv$,}\\
(iii)&\text{$K_X\sim_{\mathbb{Q}}-uf^*\mathcal{O}_{\mathbb{P}^1}(1)$}\\
\end{array}\right\}.
\]

\begin{thm}[cf.~{\cite[Theorem 1.2]{HH}}]\label{quesmain}
We fix $d\in\mathbb{Z}_{>0}$, $u\in\mathbb{Q}_{>0}$, $v\in\mathbb{Q}_{>0}$.
Then we have the following for some $r \in\mathbb{Z}_{>0}$. 
For any locally Noetherian scheme $S$ over $\mathbb{C}$, we attain a groupoid $\mathscr{M}_{d,v,u,r}(S)$ whose objects are 
 $$\left\{
 \vcenter{
 \xymatrix@C=12pt{
(\mathcal{X},\mathscr{A})\ar[rr]^-{f}\ar[dr]_{\pi_{\mathcal{X}}}&& \mathcal{C} \ar[dl]\\
&S
}
}
\;\middle|
\begin{array}{rl}
(i)&\text{$\pi_{\mathcal{X}}$ is a flat projective morphism and $\mathcal{X}$ is a scheme,}\\
(ii)&\text{$\mathscr{A}\in\mathrm{Pic}_{\mathcal{X}/S}(S)$ such that $\mathscr{A}$ is $f$-ample,}\\
(iii)&\text{$\omega_{\mathcal{X}/S}^{[r]}$ exists as a line bundle, whose restriction  to}\\
&\text{any geometric fiber $\mathcal{X}_{\bar{s}}$ over $\bar{s}\in S$ coincides with}\\
&\text{$\mathcal{O}_{\mathcal{X}_{\bar{s}}}(rK_{\mathcal{X}_{\bar{s}}})$,}\\
(iv)&\text{$\pi_{\mathcal{X}*}\omega_{\mathcal{X}/S}^{[-lr]}$ is locally free and it generates}\\
&\text{$H^0(\mathcal{X}_{s}, \mathcal{O}_{\mathcal{X}_{s}}(-lrK_{\mathcal{X}_{s}}))$ for any point $s\in S$ and any}\\
&\text{$l\in\mathbb{Z}_{>0}$,}\\
(v)&\text{$f$ is the canonical $S$-morphism to}\\
&\text{$\mathcal{C}:=\mathbf{Proj}_{S}(\oplus\pi_{\mathcal{X}*}\omega_{\mathcal{X}/S}^{[-lr]})$ and $(\mathcal{X}_{\overline{s}},\mathscr{A}_{\overline{s}})\to \mathcal{C}_{\overline{s}} \in \mathfrak{Z}_{d, v,u}$}\\
&\text{for any geometric point $\overline{s}\in S$}
\end{array}\right\}$$
and isomorphisms are $S$-isomorphisms $\alpha\colon (\mathcal{X},\mathscr{A})\to  (\mathcal{X}',\mathscr{A}')$ such that $\alpha^*\mathscr{A}'=\mathscr{A}\otimes f^*\mathscr{B}$ for some $\mathscr{B}\in\mathrm{Pic}_{\mathcal{C}/S}(S)$.

Then $\mathscr{M}_{d,v,u,r}$ is a separated Deligne-Mumford stack of finite type over $\mathbb{C}$ with a coarse moduli space $M_{d,v,u,r}$ (cf.~\cite{Ols}).

Furthermore, there exists $w>0$ such that for any geometric point $\bar{s}\in \mathscr{M}_{d,v,u,r}$, if $\mathrm{vol}(\mathscr{A}_{\bar{s}})=w$, then $\mathscr{A}_{\bar{s}}$ is ample and the object $(\mathcal{X}_{\bar{s}},\mathscr{A}_{\bar{s}})$ corresponding to $\bar{s}$ is specially $\mathrm{K}$-stable.
\end{thm}

Take $w$ as in Theorem \ref{quesmain}.
We can set the {\it CM line bundle $\Lambda_{\mathrm{CM},w}$ on $M_{d,v,u,r}$ with respect to the volume} $w$ as follows.
As \cite[Remark 6.5]{HH}, we can put the universal family $\pi_{\mathscr{U}}\colon(\mathscr{U},\mathscr{A})\to\mathscr{M}_{d,v,u,r}$ such that $\mathrm{vol}(\mathscr{A}_{\bar{s}})=w$ for any geometric fiber $(\mathscr{U}_{\bar{s}},\mathscr{A}_{\bar{s}})$ over $\mathscr{M}_{d,v,u,r}$.
Then we can define the CM line bundle $\lambda_{\mathrm{CM},\pi_{\mathscr{U}}}$ with respect to $\mathscr{A}$ on $\mathscr{M}_{d,v,u,r}$ by the construction of $(\mathscr{U},\mathscr{A})$ (cf.~\cite[Example 2.13]{HH}).
 Let $\pi\colon\mathscr{M}_{d,v,u,r}\to M_{d,v,u,r}$ be the canonical morphism to its coarse moduli space.
By \cite[Theorem 10.3]{alper}, we obtain a $\mathbb{Q}$-line bundle $\Lambda_{\mathrm{CM},w}$ on $M_{d,v,u,r}$ such that $\pi^*\Lambda_{\mathrm{CM},w}=\lambda_{\mathrm{CM},\pi_{\mathscr{U}}}.$ 
Note that any geometric fiber $(\mathscr{U}_{\bar{s}},\mathscr{A}_{\bar{s}})$ is specially K-stable by Theorem \ref{quesmain}.

\section{Nefness of the CM line bundle}

We first discuss the nonnegativity of J-degree.

\begin{prop}\label{prop--nef-j-st}
    Let $\pi\colon(X,\Delta,L)\to C$ be a polarized log family, where $C$ is a proper smooth curve, and let $H$ be a nef $\mathbb{Q}$-line bundle on $X$.
    Suppose that for any very general closed point $s\in C$, $(X_s,L_s)$ is $\mathrm{J}^{H_s}$-semistable. Then 
    $$\mathcal{J}^{H}((X,L)/C)\ge0.$$
\end{prop}

\begin{proof}
    We first deal with the case when $H$ is ample and for any very general point $s\in C$, $(X_s,L_s)$ is $\mathrm{J}^{H_s}$-semistable.
    Suppose that $rH$ is very ample.
    Take a general member $D\in|rH|$ such that $D$ is normal, $\pi|_D\colon D\to C$ is a contraction and $D_s$ is compatible with $F_{\mathrm{HN}}$, where $r\in \mathbb{Z}_{>0}$. 
    We set the Harder-Narasimhan filtrations $\mathscr{F}_{\mathrm{HN}}$ on $\oplus_{m\ge0} \pi_*\mathcal{O}_{X}(mrL)$ and $\mathscr{F}_{D,\mathrm{HN}}$ on $\oplus_{m\ge0} \pi_*\mathcal{O}_{D}(mrL|_D)$ respectively. 
    There exists the canonical map 
    \begin{equation}\label{eq---restrict--hn}
    \mathscr{F}^\lambda_{\mathrm{HN}}\pi_*\mathcal{O}_{X}(mrL)\to\mathscr{F}^\lambda_{D,\mathrm{HN}}\pi_*\mathcal{O}_{D}(mrL|_D)
\end{equation}
    for any sufficiently divisible $m\in\mathbb{Z}_{>0}$ and $\lambda\in\mathbb{Q}$ by \cite[Lemma 1.3.3]{HL}.
We consider the restricted filtration $(F_{\mathrm{HN}})|_{D_s}$ (cf.~\cite[Example 2.4]{CM}) and the induced filtration $F_{D,\mathrm{HN}}$ defined for the family $(D,rL|_D)\to S$.
We see by (\ref{eq---restrict--hn}) that for any sufficiently divisible $m\in\mathbb{Z}_{>0}$ and $\lambda\in\mathbb{Q}$, there exists a natural inclusion
\[
(F_{\mathrm{HN}})|_{D_s}^{\lambda}H^0(D_s,\mathcal{O}_{D_s}(mL_s))\subset F_{D,\mathrm{HN}}^{\lambda}H^0(D_s,\mathcal{O}_{D_s}(mL_s)).
\]
Let $w_{F_{\mathrm{HN}}}(m)$, $w_{F_{D,\mathrm{HN}}}(m)$ and $w_{(F_{\mathrm{HN}})|_{D_s}}(m)$ be the weight functions of $F_{\mathrm{HN}}$, $F_{D,\mathrm{HN}}$ and $(F_{\mathrm{HN}})|_{D_s}$ respectively.
Then we have by Proposition \ref{prop-cm-df-good} that
\begin{align*}
&\mathcal{J}^{rH}((X,L)/C)=\lim_{m\to\infty}\left(n!\frac{w_{F_{D,\mathrm{HN}}}(m)}{m^n}-\frac{nr(H_s\cdot L_s^{n-1})}{(n+1)L_s^n}\frac{(n+1)!w_{F_{\mathrm{HN}}}(m)}{m^{n+1}}\right)\\
&\ge\lim_{m\to\infty}\left(n!\frac{w_{(F_{\mathrm{HN}})|_{D_s}}(m)}{m^n}-\frac{nr(H_s\cdot L_s^{n-1})}{(n+1)L_s^n}\frac{(n+1)!w_{F_{\mathrm{HN}}}(m)}{m^{n+1}}\right)=\mathcal{J}^{rH,\mathrm{NA}}(F_{\mathrm{HN}}).
\end{align*}
By Proposition \ref{prop--filt--j}, $r\mathcal{J}^{H,\mathrm{NA}}(F_{\mathrm{HN}})=\mathcal{J}^{rH,\mathrm{NA}}(F_{\mathrm{HN}})\ge0$.

    We claim that the asserion in the general case follows from what we have shown in the previous paragraph.
    Indeed, suppose that $L+c\pi^*P$ is ample for some $c\in\mathbb{Q}_{>0}$ and closed point $P\in C$.
    Then we see that 
    \[
    \mathcal{J}^{H+\epsilon(L+c\pi^*P)}((X,L)/C)\ge0
    \]
    for any $\epsilon\in\mathbb{Q}_{>0}$.
    Thus,
     \[
    \mathcal{J}^{H}((X,L)/C)=\lim_{\epsilon\to0}\mathcal{J}^{H+\epsilon(L+c\pi^*P)}((X,L)/C)\ge0.
    \]
    We complete the proof.
\end{proof}

Next, we deal with nefness of the CM line bundle when a fiber is specially K-semistable.

\begin{prop}\label{prop--delta--nef--curve}
 Let $\pi\colon(X,\Delta,L)\to C$ be a polarized log $\mathbb{Q}$-Gorenstein family such that $C$ is a proper smooth curve and $L^{n+1}=0$.
 Take $r\in\mathbb{Z}_{>0}$ such that $rL$ is a line bundle.
 Suppose that there exist a closed point $s\in C$ and $\delta\in(0,\delta(X_s,\Delta_s,L_s)]$ such that $K_{X/C}+\Delta+\delta L$ is $\pi$-nef.
 Then the following hold.
 \begin{enumerate}
\item For any sufficiently divisible $m\in\mathbb{Z}_{>0}$ and sufficiently small $\epsilon\in\mathbb{Q}_{>0}$ and $\epsilon'\in\mathbb{Q}_{>0}$,
there exists $D\in|mrL+(m\epsilon+2g(C)) \pi^*P|$ such that $(X_s,\Delta_s+\frac{\delta+\epsilon'}{mr}D_s)$ is lc, where $g(C)$ is the genus of $C$ and $P$ is an arbitrary closed point of $C$. 
\item $K_{X/C}+\Delta+\delta L$ is nef.
 \end{enumerate}
\end{prop}

\begin{proof}
    Consider the induced filtration $F_{\mathrm{HN}}$ of $R:=\oplus_{m\ge0}R_m$, where we let $R_m=H^0(X_s,\mathcal{O}_{X_s}(mrL_s))$.
By $L^{n+1}=0$ and Proposition \ref{prop-cm-df-good}, we see that $S(F_{\mathrm{HN}})=0$.
By Theorem \ref{thm--xz--4.5} and Lemma \ref{lem--xz--4.13}, for any positive sufficiently small rational numbers $\epsilon$ and $\epsilon'$, it holds that 
$$\mathrm{lct}(X_s,\Delta_s;I_{m,-m\epsilon}(F_{\mathrm{HN}}))\ge\frac{\delta(X_s,\Delta_s,L_s)+\epsilon'}{mr}$$
for any sufficiently large $m\in\mathbb{Z}_{>0}$, where $I_{m,\lambda}(F_{\mathrm{HN}})$ is the base ideal of ${F}_{\mathrm{HN}}^\lambda R_m$.
In particular, $I_{m,-m\epsilon}(F_{\mathrm{HN}})\ne0$.
We know that $I_{m,-m\epsilon}(F_{\mathrm{HN}})$ is the image of $$\mathscr{F}_{\mathrm{HN}}^{-m\epsilon}\pi_*\mathcal{O}_X(mrL)\otimes\mathcal{O}_{X_s}(-mrL_s)\to\mathcal{O}_{X_s}.$$
By \cite[Proposition 5.7]{CP}, we see that \[
\mathscr{F}_{\mathrm{HN}}^{-m\epsilon}\pi_*\mathcal{O}_X(mrL)\subset \mathrm{Image}(H^0(X,mrL+(m\epsilon+2g(C)) \pi^*P)\otimes\mathcal{O}_C\to \pi_*\mathcal{O}_X(mrL))
\]
for any sufficiently divisible $m\in\mathbb{Z}_{>0}$.
This means that there exists an effective divisor $D\in|mrL+(m\epsilon+2g(C)) \pi^*P|$ such that the section corresponding to $D_s$ is contained in $I_{m,-m\epsilon}(F_{\mathrm{HN}})$.
Thus, we see that $(X_s,\Delta_s+\frac{\delta+\epsilon'}{mr}D_s)$ is lc.
Note that
\[
K_{X/C}+\Delta+\frac{\delta+\epsilon'}{mr}D\sim_{\mathbb{Q}}K_{X/C}+\Delta+(\delta+\epsilon')\left(L+\frac{m\epsilon+2g(C)}{mr} \pi^*P\right)
\]
and hence $K_{X/C}+\Delta+\frac{\delta+\epsilon'}{mr}D$ is nef by \cite[Theorem 1.11]{fujino-semi-positivity} (cf.~the proof of \cite[Corollary 2.14]{PX}).
Here, we take $m$ sufficiently divisible and $\epsilon,\epsilon'$ sufficiently small.
Thus, $K_{X/S}+\Delta+\delta L$ is nef by taking the limit.
We complete the proof.
\end{proof}

\begin{thm}\label{thm--cm--nefness}
Let $\pi\colon(X,\Delta,L)\to C$ be a polarized log $\mathbb{Q}$-Gorenstein family, where $C$ is a proper smooth curve.
    If there exists a closed point $s\in C$ such that $(X_s,\Delta_s,L_s)$ is specially $\mathrm{K}$-semistable and $K_{X/C}+\Delta+\delta(X_s,\Delta_s,L_s)L$ is $\pi$-nef, then 
    $$\mathrm{CM}((X,\Delta,L)/C)\ge0.$$
\end{thm}

\begin{proof}
By taking some $c\in\mathbb{Q}$ and replacing $L$ by $L+c\pi^*P$, where $P$ is a closed point of $C$, we may assume that $L^{n+1}=0$.
Then we see that 
\[
\mathrm{CM}((X,\Delta,L)/C)=\mathcal{J}^{K_{X/C}+\Delta+\delta(X_s,\Delta_s,L_s) L}((X,L)/C).
\]
Thus, the assertion follows from Propositions \ref{prop--nef-j-st} and \ref{prop--delta--nef--curve}.
\end{proof}

\section{Bigness of the CM line bundle}

First, we deal with Proposition \ref{prop--j--ample}.
To show this, we assert the following.
\begin{prop}\label{prop--j--big}
   Let $\pi\colon(X,L)\to S$ be a polarized family. 
   Suppose that $S$ is projective and there exists an ample $\mathbb{Q}$-line bundle $H$ on $X$.
   If there exists a closed point $s_0\in S$ such that $(X_{s_0},L_{s_0})$ is $\mathrm{J}^{H_{s_0}}$-semistable, then $\lambda_{\mathrm{J},\pi,H}$ is big. 
\end{prop}

\begin{proof}
Take an ample line bundle $M$ on $S$.
Take $\delta\in\mathbb{Q}_{>0}$ such that $H-\delta \pi^*M$ is also ample.
By Proposition \ref{prop--generic-j-ss}, we have that $(X_{s},L_{s})$ is $\mathrm{J}^{H_s}$-semistable for any very general closed point $s\in S$.
For any movable curve $C\to S$,
\[
\lambda_{\mathrm{J},\pi,H-\delta \pi^*M}\cdot C\ge0
\]
by Proposition \ref{prop--nef-j-st}.
By \cite{BDPP}, we conclude that $\lambda_{\mathrm{J},\pi,H-\delta \pi^*M}$ is pseudo-effective.
On the other hand, 
\[
\lambda_{\mathrm{J},\pi, \pi^*M}\cdot C=(n+1)M\cdot C
\]
for any movable curve $C\to S$.
Thus, $\lambda_{\mathrm{J},\pi,\pi^*M}$ is big by \cite{BDPP}.
Since $\lambda_{\mathrm{J},\pi,H}=\lambda_{\mathrm{J},\pi,H-\delta \pi^*M}+\delta\lambda_{\mathrm{J},\pi,\pi^*M}$, we have that $\lambda_{\mathrm{J},\pi,H}$ is big.
\end{proof}

\begin{proof}[Proof of Proposition \ref{prop--j--ample}]
    This immediately follows from Propositions \ref{prop--nef-j-st} and \ref{prop--j--big}.
    Indeed, we see that $\lambda_{\mathrm{J},\pi,H}|_V$ is big and nef for any subvariety $V\subset S$.
    Then $\lambda_{\mathrm{J},\pi,H}$ is ample by the Nakai-Moishezon criterion.  
\end{proof}

Next, we deal with Theorem \ref{thm--main--2}.
For this, we show the following technical result.

\begin{prop}\label{prop--delta-big}
    Let $\pi\colon(X,\Delta,L)\to S$ be a polarized log $\mathbb{Q}$-Gorenstein family of relative dimension $n$ with maximal variation, where $S$ is projective and $(X_{\bar{s}},\Delta_{\bar{s}})$ is klt for any geometric point $\bar{s}\in S$.
    Suppose that $\pi_*L^{n+1}\equiv0$.
    Suppose that there exists $\lambda\in\mathbb{Q}_{>0}$ such that $\lambda<\delta(X_s,\Delta_s,L_s)$ for any very general closed point $s\in S$ and $K_{X/S}+\Delta+\lambda L$ is $\pi$-ample.
Then the $\mathbb{Q}$-line bundle $\pi_*(K_{X/S}+\Delta+\lambda L)^{n+1}$ is big.
\end{prop}

\begin{proof}
We follow the argument of \cite[Lemma 7.4]{XZ}.
By taking a resolution of singularities of $S$, we may assume that $S$ is smooth.
    Take a big line bundle $H$ on $S$ and $r\in\mathbb{Z}_{>0}$ such that $M:=r(K_{X/S}+\Delta+\lambda L)$ is a line bundle.
    Due to \cite[Theorem 6.6]{XZ}, the following holds. 
    Choose a suitable $d\in\mathbb{Z}_{>0}$ and let $D=\mathrm{Supp}(\Delta)$.
    Let $W=\pi_*\mathcal{O}_X(M)$ and $Q:=\pi_*\mathcal{O}_{X}(dM)\oplus (\pi|_D)_*\mathcal{O}_{D}(dM|_D)$ and set the ranks of them $w$ and $q$ respectively.
    We note that $Q$ is not locally free in general but there exists a big open subset $S^{\circ}\subset S$ such that $D|_{S^{\circ}}$ and any irreducible component of $D|_{S^{\circ}}$ are flat over $S^{\circ}$.
    We may assume that $H^i(X_s,\mathcal{O}_{X_s}(M_s))=0$ for any $s\in S$ and $H^i(D_s,\mathcal{O}_{D_s}(M_s|_{D_s}))=0$ for any $s\in S^\circ$ and $i>0$.
    We see that $Q|_{S^{\circ}}$ is a locally free sheaf of rank $q$ and set $B$ as a Weil divisor on $S$ such that $\mathrm{det}(Q|_{S^{\circ}})\sim B|_{S^\circ}$.
    Since $S$ is smooth, we regard $B$ as a Cartier divisor.
Then, we see that there exist $m\in\mathbb{Z}_{>0}$ and a non-zero map
\[
\mathrm{Sym}^{dqm}(W^{\oplus 4w})\to\mathcal{O}_S(mB-H).
\]
For any movable curve $g\colon C\to S$, the image of $C$ contains a very general point of $S$ and hence $W$ is a nef vector bundle by Proposition \ref{prop--delta--nef--curve} and \cite[Theorem 1.11]{fujino-semi-positivity}.
This means that $g^*\mathcal{O}_S(mB-H)$ is also nef since $\mathrm{Sym}^{dqm}(W^{\oplus 4w})$ is nef (cf.~\cite[Theorem 6.1.15]{Laz2}).
Thus $B$ is big by \cite{BDPP}.

In this paragraph, we show the inequality (\ref{eq--1}) below, which is a key step to show Proposition \ref{prop--delta-big}.
Consider the following map
\[
\mathrm{det}(\pi_*\mathcal{O}_{X}(dM))\hookrightarrow\bigotimes^{q_1}_{i=1}\pi_*\mathcal{O}_X(dM)\text{ and }\mathrm{det}(\pi_*\mathcal{O}_{D}(dM|_D)|_{S^\circ})\hookrightarrow\bigotimes^{q_2}_{i=1}\pi_*\mathcal{O}_D(dM|_D)|_{S^\circ},
\]
where $q_1$ and $q_2$ are the ranks of $\pi_*\mathcal{O}_{X}(dM)$ and $\pi_*\mathcal{O}_{D}(dM|_D)$ respectively.
By them, we obtain the following embedding
\begin{equation}
    \mathrm{det}(Q|_{S^\circ})\hookrightarrow\bigotimes^{q_1}_{i=1}\pi_*\mathcal{O}_X(dM)|_{S^\circ}\otimes\bigotimes^{q_2}_{i=1}\pi_*\mathcal{O}_D(dM|_D)|_{S^\circ}.\label{eq--adj-before}
\end{equation}
Let $Z:=X^{(q_1)}\times_SD^{(q_2)}$, where $X^{(q_1)}:=X\times_SX\times_S\ldots\times_SX$ means the $q_1$-times self fiber product of $X$ over $S$.
Let $M_Z:=\sum_{i=1}^{q_1}p_i^*M+\sum_{j=1}^{q_2}p'^*_jM|_D$, where $p_i\colon Z\to X$ is the $i$-th projection and $p'_j\colon Z\to D$ is the $q_1+j$-th projection. 
Let $f\colon Z\to S$ denote the canonical morphism.
Then we see that (see \cite[\S2.2]{CP}) 
\[
\bigotimes^{q_1}_{i=1}\pi_*\mathcal{O}_X(dM)\otimes\bigotimes^{q_2}_{i=1}\pi_*\mathcal{O}_D(dM|_D)\cong f_*\mathcal{O}_Z(dM_Z).
\]
By the adjunction of $f_*$ and $f^*$ applied to (\ref{eq--adj-before}), we have a non-zero map
\[
f^*\mathcal{O}_S(B)|_{f^{-1}(S^\circ)}\to\mathcal{O}_{Z}(dM_Z)|_{f^{-1}(S^\circ)}.
\]
This means that $(dM_Z-f^*B)|_{f^{-1}(S^\circ)}$ is effective on some irreducible component of $f^{-1}(S^\circ)$.
Recall that any irreducible component of $D\cap\pi^{-1}(S^\circ)$ is flat over $S^{\circ}$.
Thus, so is $f^{-1}(S^\circ)$ and hence we see that any irreducible component of $f^{-1}(S^\circ)$ can be denoted as $\pi^{-1}(S^\circ)^{(q_1)}\times_{S^\circ}\pi^{-1}(S^\circ)\cap D_1\times_{S^\circ}\ldots\times_{S^\circ}\pi^{-1}(S^\circ)\cap D_{q_2}$ for some irreducible components $D_1,\ldots, D_{q_2}$ of $D$.
We can also check that $f^{-1}(S^\circ)$ is generically reduced.
Let $Z':=X^{(q_1)}\times_SD_1\times_S\ldots\times_SD_{q_2}$ and $Z'_1$ the Zariski closure in $Z'$ $$\overline{\pi^{-1}(S^\circ)^{(q_1)}\times_{S^\circ}\pi^{-1}(S^\circ)\cap D_1\times_{S^\circ}\ldots\times_{S^\circ}\pi^{-1}(S^\circ)\cap D_{q_2}}.$$
Let $\iota\colon Z'_1\hookrightarrow Z'$ be the natural inclusion and $\nu\colon Z'^{\nu}\to Z_1'$ the normalization.
Since $\mathrm{codim}_Z(Z\setminus{f^{-1}(S^\circ)})\ge2$, $\mathrm{codim}_{Z'^\nu}(Z'^\nu\setminus{\nu^{-1}(f^{-1}(S^\circ)}\cap Z'_1))\ge2$. 
By the $S_2$-condition of $Z'^\nu$, there exists a non-zero map
\begin{equation}
    \nu^*\iota^*f|_{Z'}^*\mathcal{O}_S(B)\to\nu^*\iota^*\mathcal{O}_{Z'_1}(dM_Z|_{Z'}).\label{eq--psef--engine}
\end{equation}
We denote the base change by $g$ for any movable curve $g\colon C\to S$ of $f|_{Z'}$, $M_Z|_{Z'}$, $\iota\colon Z'_1\hookrightarrow Z'$ and $\nu\colon Z'^\nu\to Z'_1$ by $f_{Z'_C}$, $M_{Z'_C}$, $\iota_{C}$ and $\nu_C$.
Let $B_C:=g^*B$.
We note that $M_{Z'_C}$ is nef by Proposition \ref{prop--delta--nef--curve}.
By the property of (\ref{eq--psef--engine}),  $d\nu_C^*\iota_C^*M_{Z'_C}-\nu_C^*\iota_C^*f_{Z'_C}^*B_C$ is effective 
for any movable curve $C\to S$.
Thus, we obtain that \begin{equation}(d\nu_C^*\iota_C^*M_{Z'_C}-\nu_C^*\iota_C^*f_{Z'_C}^*B_C)\cdot\nu_C^* \iota_C^*M_{Z'_C}^{N-1}\ge0.\label{eq--addition}
\end{equation}
Let $N=dq_1+(d-1)q_2=\mathrm{dim}\,Z'_{1,C}-1$.
Then we have that $\mathrm{dim}\,(Z'_C\setminus Z'_{1,C})\le N$ since each fiber of each $D_i\to S$ is of dimension at most $n-1$.
This means that for any $N+1$ line bundles $L_1,L_2,\ldots,L_{N+1}$ on $Z'_C$, 
\[
L_1\cdot\ldots\cdot L_{N+1}=L_{1}|_{Z'_{1,C}}\cdot\ldots\cdot L_{N+1}|_{Z'_{1,C}}.
\]
Therefore, we have by (\ref{eq--addition}) that
\begin{equation}
dM_{Z'_C}^{N+1}\ge (M_{Z'_C,t}^N)\mathrm{deg}_C\,B_C.\label{eq--1}
\end{equation}

To complete the proof of Proposition \ref{prop--delta-big}, we have to show by (\ref{eq--1}) that there exists a positive constant $C_4>0$ such that $M_{C}^{n+1}\ge C_4\mathrm{deg}_C\,B_C$ for any movable curve $C\to S$.
Let $$C_0:=\max\left\{(M_t^n)^{q_1}\prod_{i=1}^{q_2}(M_t|_{D'_{i,t}})^{n-1}\right\}>0$$ be a constant, where $D'_1,\ldots ,D'_{q_2}$ run over all $q_2$ irreducible components of $D|_{\pi^{-1}(S^\circ)}$.
Here, we note that $(M_{Z_C',t}^{N})=(M_t^n)^{q_1}\prod_{i=1}^{q_2}(M_t|_{D_{i,t}})^{n-1}$ and thus $(M_{Z_C',t}^{N})\ge C_0$.
Next, we see as \cite[(6.3.5.i)]{P} that there exists a constant $C_1>0$ such that 
\begin{equation}
    M_C^{n+1}+M_{C}^n\cdot\Delta_C\ge C_1M_{Z'_C}^{N+1} \label{eq--2}
\end{equation}
independent from the choice of $D_1,\ldots,D_{q_2}$ and $C$.
Indeed, let $D_{i,C}:=D_i\times_SC$ and take the Zariski closure $D^*_{i,C}:=\overline{D_{i,C}\cap \pi^{-1}(S^\circ)}\subset D_{i,C}$ for each $1\le i\le q_2$.
It is easy to see that $D^*_{i,C}$ is flat over $C$.
Let  
$$Z_2':=X_C^{(q_1)}\times_CD^*_{1,C}\times_C\ldots\times_CD^*_{q_2,C}\subset Z'_C.$$
Since each fiber of each $D_i\to S$ is of dimension at most $n-1$, $(M_{Z'_C})^{N+1}=(M_{Z'_C}|_{Z'_2})^{N+1}$. 
By applying \cite[Lemma 7.0.5]{P} to $Z'_2$, there exists a positive constant $d_1\in\mathbb{Q}_{>0}$ depending only on $n$ such that 
\begin{align*}(M_{Z'_C})^{N+1}&=q_1d_1(M_C^{n+1})(M_t^n)^{q_1-1}\prod_{i=1}^{q_2}(M_t|_{D_{i,t}})^{n-1}\\
&+ d_1\sum_{i=1}^{q_2}(M_{C}|_{D_i})^n(M_t^n)^{q_1}\prod_{j\ne i}^{q_2}(M_t|_{D_{j,t}})^{n-1}.
\end{align*}
This proves the existence of such $C_1$.
Thus, we see by (\ref{eq--1}) and (\ref{eq--2}) that there exists a positive constant $C_2:=d^{-1}C_0C_1$ independent from the choice of movable curves $C$ such that
\begin{equation}
  M_C^{n+1}+M_{C}^n\cdot\Delta_C\ge C_2\mathrm{deg}_C\,B_C.\label{eq--3}  
\end{equation}
Now, it suffices to show the following claim.
\begin{claim}\label{cl--xz}
There exists a positive constant $C_3>0$ independent from the choice of movable curves $C\to S$ such that
\[
M_C^{n+1}\ge C_3(M_{C}^n\cdot\Delta_C).
\]
\end{claim}
\begin{proof}[Proof of Claim \ref{cl--xz}]
We mimic the proof of \cite[Lemma 7.6]{XZ}.
We note that there exists $0<\xi<1$ such that $K_{X_s}+(1-\xi)\Delta_s+\lambda L_s$ is big as a $\mathbb{Q}$-Weil divisor for any very general point $s\in S$.
Indeed, we choose $\xi$ such that $K_{X_\eta}+(1-\xi)\Delta_\eta+\lambda L_\eta$ is big, where $\eta$ is the generic point of $S$.
Then we see that $K_{X_s}+(1-\xi)\Delta_s+\lambda L_s$ is big for general $s$. 
For any movable curve $C\to S$, we see that $(X_C,\Delta_C)$ is klt since $(X_c,\Delta_c)$ is klt for any closed point $c\in C$.
Thus, there exists a small projective birational morphism $g\colon Y\to X_C$ from a normal $\mathbb{Q}$-factorial variety by \cite[1.4.3]{BCHM}.
Let $\Delta'_C:=g_*^{-1}\Delta_C$ and let $\varphi\colon Y\to C$ be the canonical morphism.
Fix $r'\in\mathbb{Z}_{>0}$ such that $r'L$ is a line bundle.
By Proposition \ref{prop--delta--nef--curve}, we see that for any sufficiently small $\epsilon>0$ and sufficiently divisible $m\in\mathbb{Z}_{>0}$,
there exists an effective divisor $D\in|mr'L_C+(m\epsilon+2g(C)) f^*P|$ such that $\mathrm{lct}(X_s,\Delta_s;D_s)\ge\frac{\lambda+\epsilon}{mr'}$, where $g(C)$ is the genus of $C$, $P\in C$ is a closed point and $s\in C$ is a very general point.
Thus, we see that $(X_s,\Delta_s+\frac{\lambda+\epsilon}{mr'}D_s)$ is lc for any sufficiently general $s\in C$. 
Let $\Gamma:=(1-\xi)\Delta'_C+\frac{\lambda+\epsilon}{mr'}g^*D$.
 Then  $K_{Y_s}+\Gamma_s$ is big for any sufficiently general $s\in C$ and hence \begin{equation}
     \label{eq--big--fiberwise}H^0(Y_s,\mathcal{O}_{Y_s}(l(K_{Y_s}+\Gamma_s)))\ne0
 \end{equation} 
 for any sufficiently divisible $l\in\mathbb{Z}_{>0}$.
 Let $\psi\colon Y_{\mathrm{lc}}\to Y$ be the lc modification of $(Y,\Gamma)$ by \cite[Theorem 1.1]{OX}.
In other words, $\psi$ is a projective birational morphism of normal varieties and  there exists an effective $\psi$-exceptional $\mathbb{Q}$-divisor $G$ such that 
\[
\psi^*(K_Y+\Gamma)-G=K_{Y_{\mathrm{lc}}}+\psi_*^{-1}\Gamma+\mathrm{Ex}(\psi),
\]
  $K_{Y_{\mathrm{lc}}}+\psi_*^{-1}\Gamma+\mathrm{Ex}(\psi)$ is $\psi$-ample and $(Y_{\mathrm{lc}},\psi_*^{-1}\Gamma+\mathrm{Ex}(\psi))$ is lc.
  Since $(X_s,\Delta_s+\frac{\lambda+\epsilon}{mr'}D_s)$ is lc for sufficiently general $s\in C$, we have that $G$ is vertical with respect to $C$.
  Therefore, there exists a coherent sheaf $\mathcal{G}_l$ on $C$ whose support is zero-dimensional for any sufficiently divisible $l\in\mathbb{Z}_{>0}$ such that there exists an exact sequence
  \[
  0\rightarrow (\varphi\circ\psi)_*\mathcal{O}_{Y_{\mathrm{lc}}}(l(K_{Y_{\mathrm{lc}}/C}+\psi_*^{-1}\Gamma+\mathrm{Ex}(\psi)))\rightarrow \varphi_*\mathcal{O}_Y(l(K_{Y/C}+\Gamma))\rightarrow\mathcal{G}_l\rightarrow0.
    \]
    By \cite[Theorem 1.1]{fujino--semipositive--another}, we have that $(\varphi\circ\psi)_*\mathcal{O}_{Y_{\mathrm{lc}}}(l(K_{Y_{\mathrm{lc}}/C}+\psi_*^{-1}\Gamma+\mathrm{Ex}(\psi)))$ is weakly positive over $C$ for any sufficiently divisible $l$.
    Since $\mathrm{dim}\,\mathrm{Supp}\,\mathcal{G}_l=0$, $\varphi_*\mathcal{O}_Y(l(K_{Y/C}+\Gamma))$ is also weakly positive.
    This means that for any ample line bundle $A$ on $C$ and positive integer $a$, there exists $b\in\mathbb{Z}_{>0}$ such that the stalk of $\mathrm{Sym}^{ab}(g_*\mathcal{O}_Y(l(K_{Y/C}+\Gamma)))\otimes\mathcal{O}_C(bA)$ at the generic point of $C$ is generated by $H^0(C,\mathrm{Sym}^{ab}(g_*\mathcal{O}_Y(l(K_{Y/C}+\Gamma)))\otimes\mathcal{O}_C(bA))$. 
    By the following commutative diagram
    \[
    \begin{CD}
\mathrm{Sym}^{ab}(\varphi_*\mathcal{O}_Y(l(K_{Y/C}+\Gamma)))\otimes\mathcal{O}_C(bA)@>>>\varphi_*\mathcal{O}_Y(abl(K_{Y/C}+\Gamma))\otimes\mathcal{O}_C(bA)\\
@VVV @VVV\\
\mathrm{Sym}^{ab}H^0(Y_s,\mathcal{O}_{Y_s}(l(K_{Y_s}+\Gamma_s)))@>>> H^0(Y_s,\mathcal{O}_{Y_s}(abl(K_{Y_s}+\Gamma_s))),
\end{CD}
    \]
     (\ref{eq--big--fiberwise})
    and the fact that \begin{align*}
        \mathrm{Sym}^{ab}(\varphi_*\mathcal{O}_Y(l(K_{Y/C}+\Gamma)))\otimes(\mathcal{O}_C/\mathfrak{m}_s)&\cong\mathrm{Sym}^{ab}H^0(Y_s,\mathcal{O}_{Y_s}(l(K_{Y_s}+\Gamma_s)))\\
        \varphi_*\mathcal{O}_Y(abl(K_{Y/C}+\Gamma))\otimes(\mathcal{O}_C/\mathfrak{m}_s)&\cong H^0(Y_s,\mathcal{O}_{Y_s}(abl(K_{Y_s}+\Gamma_s))),
    \end{align*} 
    where $\mathfrak{m}_s$ is the maximal ideal sheaf corresponding to $s$, for any very general point $s\in C$,
 $H^0(C,g_*\mathcal{O}_Y(abl(K_{Y/C}+\Gamma))\otimes\mathcal{O}_C(bA))\ne0$.
    This means that $al(K_{Y/C}+\Gamma)+g^*A$ is effective.
    By considering $a\to\infty$, we obtain that $K_{Y/C}+\Gamma$ is pseudo-effective and hence so is $K_{X_C/C}+(1-\xi)\Delta_C+\frac{\lambda+\epsilon}{mr'}D$.
   Therefore, we obtain that $K_{X_C/C}+(1-\xi)\Delta_C+\lambda L_C$ is pseudo-effective.
   This means that 
   \[
M_C^{n+1}\ge  r\xi (M_{C}^n\cdot \Delta_C).
\]
By taking $C_3=r\xi$, we complete the proof of Claim \ref{cl--xz}.
    \end{proof}
    By (\ref{eq--3}), Claim \ref{cl--xz} and \cite{BDPP}, we obtain that there exists a positive constant $C_4$ such that $\pi_*(M^{n+1})-C_4B$ is pseudo-effective.
    Since $B$ is big, so is $\pi_*(M^{n+1})$.
    We complete the proof.
    \end{proof}
    
\begin{proof}[Proof of Theorem \ref{thm--delta--ampleness}]
    This immediately follows from Propositions \ref{prop--delta--nef--curve} and \ref{prop--delta-big} by the Nakai--Moishezon criterion.
\end{proof}    
By applying Proposition \ref{prop--delta-big}, we show the following key ingredient to prove Theorem \ref{thm--main--1}.

\begin{thm}\label{thm--CM-bigness}
    Let $\pi\colon(X,\Delta,L)\to S$ be a polarized log $\mathbb{Q}$-Gorenstein family with maximal variation, where $S$ is projective and $(X_{\bar{s}},\Delta_{\bar{s}})$ is klt for any geometric point $\bar{s}\in S$.
    Suppose that there exists a closed point $s_0\in S$ such that $(X_{s_0},\Delta_{s_0},L_{s_0})$ is specially {\rm K}-stable and $K_{X/S}+\Delta+\delta(X_{s_0},\Delta_{s_0},L_{s_0})L$ is $\pi$-ample.
Then the CM-line bundle $\lambda_{\mathrm{CM},\pi}$ is big.
\end{thm}

\begin{proof}
    Let $n$ be the relative dimension of $\pi$ and $v=L_{s_0}^n$.
    Then, for any movable curve $C\to S$, the pullback of $(n+1)vL-\pi^*(\pi_*(L^{n+1}))$ satisfies 
    \[
    \left((n+1)vL-\pi^*(\pi_*(L^{n+1})\right)_C^{n+1}=0.
    \]
    Thus, we may assume that $L_C^{n+1}=0$ for any movable curve $C\to S$ by replacing $L$ with  $(n+1)vL-\pi^*(\pi_*(L^{n+1}))$.
    Next, we take positive rational numbers $\lambda$ and $\epsilon$ such that for any very general point $s\in S$, $\delta_{(X_s,\Delta_s,L_s)}\ge\lambda +\epsilon$ and $(X_s,L_s)$ is $\mathrm{J}^{K_{X_s}+\Delta_s+\lambda L_s}$-semistable by Corollary \ref{cor--verygeneral-special-kst}. 
    Here, we may assume that $K_{X/S}+\Delta+\lambda L$ is $\pi$-ample.
    By taking a suitable $r\in\mathbb{Z}_{>0}$, we may further assume that $M:=r(K_{X/S}+\Delta+\lambda L)$ is a $\pi$-very ample line bundle.

By Proposition \ref{prop--delta-big}, $\pi_*(M^{n+1})$ is big.
This means that for any movable curve $C\to S$, $M_C^{n+1}>0$.
If we choose $0<\delta<\frac{1}{(n+1)(M_t^n)}$, then we see by \cite[Theorem 2.2.15]{Laz} that $$M_{C}-\delta(\pi_C)^*(\pi_C)_*(M_{C}^{n+1})=r\left(K_{X_C/C}+\Delta_C+\lambda L_C-\frac{\delta}{r}\pi_C^*(\pi_C)_*(M_{C}^{n+1})\right)$$ is big, where $t$ is a general closed point of $S$.
Here, we claim the following.
\begin{claim}\label{claim--alpha}
Let $\alpha:=\inf_{t\in S}\alpha(X_t,\Delta_t;M_t)$.
Then $\alpha>0$ and
\begin{align}
   K_{X_C/C}+\Delta_C+\lambda L_C-\frac{\alpha\delta\epsilon}{\lambda+(1+r\alpha)\epsilon}\pi_C^*(\pi_C)_*(M_C^{n+1})\label{eq-claim-2}
\end{align} 
is nef.
\end{claim}
\begin{proof}[Proof of Claim \ref{claim--alpha}]
Take $D'\in|M_{C}-\delta\pi_C^*(\pi_C)_*(M_{C}^{n+1})|_{\mathbb{Q}}$ by the bigness and assume that $\mathrm{Supp}\,D'$ does not contain $X_s$ for some very general closed point $s\in C$.
By Proposition \ref{prop--delta--nef--curve}, we see that for any sufficiently small $\eta\in\mathbb{Q}_{>0}$ and sufficiently divisible $m\in\mathbb{Z}_{>0}$,
there exists an effective divisor $D\in|mL_C+(m\eta+2g(C)) \pi_C^*P|$ such that $\mathrm{lct}(X_s,\Delta_s;D_s)\ge\frac{\lambda+\epsilon}{m}$, where $g(C)$ is the genus of $C$, $P\in C$ is a closed point and $s\in C$ is a very general point.
Then, we have that for any prime divisor $E$ over $X_s$,
\begin{align*}
    \frac{\lambda}{(\lambda+\epsilon)}A_{(X_s,\Delta_s)}(E)\ge\frac{\lambda}{m}\mathrm{ord}_E(D_s). 
\end{align*}
On the other hand, $\alpha>0$ by \cite[Proposition 5.3]{BL}.
Thus, we have that $(X_s,\Delta_s+\frac{\lambda}{m}D_s+\frac{\alpha\epsilon}{\lambda+\epsilon}D'_s)$ is lc.
This means that
\begin{align*}
K_{X_C/C}+\Delta_C+\frac{\lambda}{m}D+\frac{\alpha\epsilon}{\lambda+\epsilon}D'&\sim_{\mathbb{Q}}K_{X_C/C}+\Delta_C+\lambda\left(L_C+\frac{m\eta+2g(C)}{m} \pi_C^*P\right)\\
&+\frac{r\alpha\epsilon}{\lambda+\epsilon}\left(K_{X_C/C}+\Delta_C+\lambda L_C-\frac{\delta}{r}\pi_C^*(\pi_C)_*(M_{C}^{n+1})\right).
\end{align*}
is nef by \cite[Theorem 1.11]{fujino-semi-positivity}.
Since this holds for any sufficiently small $\eta$ and large $m$, we have that
(\ref{eq-claim-2})
is nef.
\end{proof}
Thus, we have that for any movable curve $C\to S$,
\begin{align*}
\mathrm{CM}((X_C,\Delta_C,L_C)/C)&=\mathcal{J}^{\left(K_{X_C/C}+\Delta_C+\lambda L_C-\frac{\alpha\delta\epsilon}{\lambda+(1+r\alpha)\epsilon}\pi_C^*(\pi_C)_*(M_C^{n+1})\right)}((X_C,L_C)/C)\\
&+\frac{\alpha\delta\epsilon}{\lambda+(1+r\alpha)\epsilon}(M_C^{n+1})v.
\end{align*}
Since $(X_s,L_s)$ is $\mathrm{J}^{K_{X_s}+\Delta_s+\lambda L_s}$-semistable for any very general point $s\in S$, we have by Proposition \ref{prop--nef-j-st} that
\[
\frac{1}{(n+1)v}\lambda_{\mathrm{CM},\pi}\cdot C\ge \frac{\alpha\delta\epsilon v}{\lambda+(1+r\alpha)\epsilon}(\pi_*(M^{n+1})\cdot C).
\]
Thus, $\lambda_{\mathrm{CM},\pi}$ is big by \cite{BDPP} since $\pi_*(M^{n+1})$ is big.
\end{proof}

\begin{proof}[Proof of Theorem \ref{thm--main--2}]
    This immediately follows from Theorems \ref{thm--cm--nefness} and \ref{thm--CM-bigness} in the same way as Proposition \ref{prop--j--ample}. 
\end{proof}

\section{An application to the moduli of K-stable Calabi-Yau fibrations over curves}
In this section, we show the following stronger result than Corollary \ref{thm--main--1}.
\begin{thm}\label{thm--proj-our-moduli}
There exists $w\in\mathbb{Z}_{>0}$ such that for any proper subspace $B$ of $M_{d,v,u,r}$, $\Lambda_{\mathrm{CM},w}|_B$ is ample.
In particular, $B$ is projective.   
\end{thm}

First, we recall the following well-known result.

\begin{lem}[cf.~{\cite[Prop.~8.3]{kawamata}}, {\cite[Prop.~4.2]{DG}}]\label{lem-pic-aut}
    Let $(X,\Delta)$ be a projective klt pair such that $K_X+\Delta\sim_{\mathbb{Q}}0$.
    Then $\mathrm{dim}\,\mathrm{Aut}_0(X,\Delta)=\mathrm{dim}\,\mathrm{Pic}^0(X)$  and for any two ample line bundles $A_1$ and $A_2$ algebraically equivalent to each other, there exists $\xi\in\mathrm{Aut}_0(X,\Delta)$ such that $\xi^*A_1\sim A_2$.
\end{lem}

\begin{proof}
For the reader's convenience, we show this lemma here.
First, we show that $\mathrm{dim}\,\mathrm{Aut}_0(X,\Delta)\le\mathrm{dim}\,\mathrm{Pic}^0(X)$.
Fix a very ample line bundle $L$ on $X$.
Consider a morphism
\[
\varphi_L\colon\mathrm{Aut}_0(X,\Delta)\ni g\mapsto [g^*L\otimes L^{\otimes-1}]\in\mathrm{Pic}^0(X).
\]
By \cite[\S4, Cor.~1]{Ab} and \cite[Prop.~4.6]{Am}, $\varphi_L$ is a homomorphism of Abelian varieties.
Thus, it suffices to show that $\mathrm{Ker}\,\varphi_L$ is a finite group scheme.
Let $\iota\colon X\hookrightarrow\mathbb{P}^{h^0(X,\mathcal{O}_X(L))-1}$ be the natural embedding defined by $|L|$.
Since $g\in \mathrm{Ker}\,\varphi_L$ satisfies that $g^*L\sim L$, there exists a group homomorphism $\nu\colon(\mathrm{Ker}\,\varphi_L)^0\to PGL(h^0(X,\mathcal{O}_X(L)))$ such that $(\mathrm{Ker}\,\varphi_L)^0$ acts on $\mathbb{P}^{h^0(X,\mathcal{O}_X(L))-1}$ so that $\iota$ is $(\mathrm{Ker}\,\varphi_L)^0$-equivariant, where $(\mathrm{Ker}\,\varphi_L)^0$ is the identity component of $\mathrm{Ker}\,\varphi_L$.
It is easy to see that $\nu$ is trivial and $(\mathrm{Ker}\,\varphi_L)^0$ trivially acts on $(X,\Delta)$.
Therefore, $\mathrm{Ker}\,\varphi_L$ is a finite group scheme.
We note that if $\mathrm{dim}\,\mathrm{Aut}_0(X,\Delta)\ge\mathrm{dim}\,\mathrm{Pic}^0(X)$, then $\varphi_L$ is further \'{e}tale. 

    We prove $\mathrm{dim}\,\mathrm{Aut}_0(X,\Delta)\ge\mathrm{dim}\,\mathrm{Pic}^0(X)$ by induction on $\mathrm{dim}\,X=n$.
    It is well-known that the assertion holds when $n=1$. 
    We may assume that $n>1$. 
    Since $(X,\Delta)$ is klt, $X$ has only rational singularities by \cite[Theorem 5.22]{KM}.
    Thus, $\mathrm{dim}\,\mathrm{Pic}^0(X)=\mathrm{dim}\,\mathrm{Alb}(X)$, where $\pi\colon X\to \mathrm{Alb}(X)$ is the Albanese morphism (cf.~\cite[\S8]{kawamata}).
    By \cite[Theorem 4.8]{Am}, we have that there exist an \'{e}tale morphism $A\to \mathrm{Alb}(X)$ from an Abelian variety, a projective connected klt log pair $(F,\Delta_F)$ and an isomorphism over $A$
    \[
    \Phi\colon A\times_{\mathrm{Alb}(X)}(X,\Delta)\to A\times (F,\Delta_F).
    \]
    Note that $A\to \mathrm{Alb}(X)$ is an \'{e}tale Galois covering and let $G=\mathrm{Ker}\,(A\to \mathrm{Alb}(X))$ be the Galois group.
    We see that $G$ is a finite commutative group.
    By identifying $(F,\Delta_F)$ with the fiber of $\pi\colon(X,\Delta)\to\mathrm{Alb}(X) $ over $0$, $G$ acts on $(F,\Delta_F)$ naturally.
    Let $\psi\colon G\to \mathrm{Aut}(F,\Delta_F)$ be the natural homomorphism induced by the $G$-action.
    On the other hand, $G$ naturally acts on $(X,\Delta)\times_{\mathrm{Alb}(X)}A$ equivariantly over $A$.
    By $\Phi$, we obtain the induced $G$-action on $A\times(F,\Delta_F)$ such that
    \[
    g\cdot(a,f)=(a+g,\phi_g(a)(f)), 
    \]
    where $g\in G$, $a\in A$ and $f\in F$ are closed points.
Here, $\phi_g(a)\in \mathrm{Aut}(F,\Delta_F)$.
Note that $\phi_g(0)=\Phi(g,\cdot)\circ\psi(g)\circ\Phi(0,\cdot)^{-1}$.
Thus, $\phi_g(0)$ is contained in the same component of $\mathrm{Aut}(F,\Delta_F)$ as $\psi(g)$.
Since $\phi_g(a)$ is continuous on $a\in A$, we can write
\[
\phi_g(a)=\psi(g)\circ t_g(a),
\]
where $t_g\colon A\to \mathrm{Aut}_0(F,\Delta_F)$ is a morphism of Abelian varieties.

If $\mathrm{Alb}(X)$ is a point, then $\mathrm{dim}\,\mathrm{Aut}_0(X,\Delta)=0$ also holds by what we have shown in the first paragraph.
   Thus, we may assume that $\mathrm{dim}\,\mathrm{Alb}(X)>0$ and then $\mathrm{dim}\,F<n$.
   Take a very ample line bundle $L$ on $X$.
   Let $\tilde{L}$ be the pullback of $L$ to $A\times F$ under the morphism $A\times F\to X$.
   $L_a$ denotes the restriction of $\tilde{L}$ to $\{a\}\times F\subset A\times F$ for any closed point $a\in A$.
   For any closed point $a\in A$ and $g\in G$, we have that 
   \begin{align}
       L_a&=g^*\tilde{L}\otimes\mathcal{O}_{\{a\}\times F}=\phi_g(a)^*L_{g+a}\label{eq--cocycle--cond--pic}\\
       &=t_g(a)^*(\psi(g)^*L_{g+a}).\nonumber
   \end{align}
   Set $$\rho\colon A\ni a\mapsto [L_a\otimes L_0^{\otimes-1}]\in \mathrm{Pic}^0(F).$$ 
   $\rho$ is indeed a morphism.
   We also consider the following morphism
   \[
   \varphi:=\varphi_{L_0}\colon\mathrm{Aut}_0(F,\Delta_F)\ni g\mapsto [g^*L_0\otimes L_0^{\otimes-1}]\in\mathrm{Pic}^0(F).
   \]
   By what we have shown in the first paragraph and the induction hypothesis, $\varphi$ is an \'{e}tale homomorphism.
   Furthermore, by \cite[Lemma 1.6]{KM}, we see that
   \[
   \varphi(h)=[h^*L_b\otimes L_b^{\otimes-1}]
   \]
   for any $h\in \mathrm{Aut}_0(F,\Delta_F)$ and $b\in A$.
   Thus (\ref{eq--cocycle--cond--pic}) is rephrasable as 
   \[
  \rho(a)-\varphi(t_g(a))=[\psi(g)^*L_{g+a}\otimes L_0^{\otimes-1}].
   \]
   Since $L_0$ is $G$-invariant, we see that
   \[
   \psi(g)^*L_{g+a}\otimes L_0^{\otimes-1}=\psi(g)^*(L_{g+a}\otimes L_0^{\otimes-1}).
   \]
   Thus we obtain that
   \begin{equation}
       \rho(a)-\varphi(t_g(a))=\psi(g)^*\rho(g+a)\label{eq--cocycle--cond-pic-2}.
   \end{equation}
   
Then, consider the following cartesian diagram
\[
 \begin{CD}
\tilde{A}_1@>>>\mathrm{Aut}_0(F,\Delta_F) \\
@VVV @V{\varphi}VV  \\
A@>_{\rho}>>\mathrm{Pic}^0(F) 
\end{CD}
\]
and let $\tilde{A}$ be the identity component of $\tilde{A}_1$.
We see that $\tilde{A}$ is an Abelian vaiety since this is a projective algebraic group.
Let $\eta\colon \tilde{A}\to A$ be the natural morphism.
Then there exists a morphism $\tilde{\rho}\colon \tilde{A}\to \mathrm{Aut}_{0}(F,\Delta_F)$ such that $\varphi\circ\tilde{\rho}=\rho\circ\eta$.
Let $H$ be a Galois group of $p\colon\tilde{A}\to\mathrm{Alb}(X)$ and let $q\colon H\to G$ be the natural morphism.
Via $q$, $H$ acts on $\tilde{A}\times (F,\Delta_F)$ equivariantly over $A\times (F,\Delta_F)$.
We denote the automorphism of $\tilde{A}\times (F,\Delta_F)$ by $b_h$ induced by $h\in H$. 
Note that $G$ acts on $\mathrm{Aut}_0(F,\Delta_F)$ and $\mathrm{Pic}^0(F)$ in the way that 
\[
g\cdot s=\psi(g)\circ s\circ\psi(g^{-1}),
\]
for any $g\in G$ and $s\in \mathrm{Aut}_0(F,\Delta_F)$, 
and 
\[
g\cdot [M\otimes L_0^{\otimes-1}]=[\psi(g^{-1})^*M\otimes L_0^{\otimes-1}]
\]
for any $g\in G$ and $[M]\in \mathrm{Pic}^0(F)$ respectively.
We see that $\varphi$ is $G$-equivariant.
Let $\tilde{t}_h(\tilde{a}):=t_{q(h)}(p(\tilde{a}))$ for any $h\in H$ and $\tilde{a}\in \tilde{A}$.
By (\ref{eq--cocycle--cond-pic-2}), if we put
\[
\theta_h(\tilde{a}):=\tilde{\rho}(\tilde{a})-\tilde{t}_h(\tilde{a})-\psi(q(h^{-1}))\circ\tilde{\rho}(\tilde{a}+h)\circ\psi(q(h))
\]
for $h\in H$ and $\tilde{a}\in \tilde{A}$, then we have that $\theta_h(\tilde{a})\in\mathrm{Ker}\,q=\mathrm{Ker}\,p$.
Since $\theta_h\colon \tilde{A}\to \mathrm{Ker}\,p$ is a morphism and $\mathrm{Ker}\,p$ is  finite, $\theta_h(\tilde{a})$ is independent of $\tilde{a}$ and we also denote $\theta_h=\theta_h(\tilde{a})\in\mathrm{Ker}\,p$.
Put an automorphism of $\tilde{A}\times (F,\Delta_F)$ over $\mathrm{Alb}(X)$ as
\[
\Psi\colon(\tilde{a},f)\mapsto (\tilde{a},\tilde{\rho}(\tilde{a})(f))
\]
and a morphism $c_h$ for any $h\in H$ as
\[
c_h\colon\tilde{A}\times F\ni(\tilde{a},f)\mapsto (\tilde{a}+h,\psi(q(h))\circ\theta_h^{-1}(f))\in\tilde{A}\times F.
\]
Then, we see that $c_h=\Psi\circ b_h\circ\Psi^{-1}$.
By $\Psi$, we may assume that $H$ acts on $\tilde{A}\times (F,\Delta_F)$ by $c_h$ and then we see that the automorphism of $\tilde{A}\times(F,\Delta_F)$
\[
\sigma_{\tilde{b}}\colon(\tilde{a},f)\mapsto (\tilde{a}+\tilde{b},f)
\]
is $H$-invariant for any $\tilde{b}\in\tilde{A}$.
This means that $\sigma_{\tilde{b}}$ descends to an automorphism of $(X,\Delta)$ and hence $\tilde{A}$ acts on $\mathrm{Alb}(X)$ transitively.
Therefore, $\mathrm{Aut}_0(X,\Delta)$ acts on $\mathrm{Alb}(X)$ transitively.
This shows that $\mathrm{dim}\,\mathrm{Pic}^0(X)=\mathrm{dim}\,\mathrm{Alb}(X)\le\mathrm{dim}\,\mathrm{Aut}_0(X,\Delta)$.
We complete the proof of $\mathrm{dim}\,\mathrm{Aut}_0(X,\Delta)=\mathrm{dim}\,\mathrm{Pic}^0(X)$.

Finally, we deal with the last assertion. 
Take $m\in\mathbb{Z}_{>0}$ such that $A_1^{\otimes m}$ is very ample.
We see that 
\[
\varphi_{A_1}\colon\mathrm{Aut}_0(X,\Delta)\ni g\mapsto [g^*A_1\otimes A_1^{\otimes-1}]\in\mathrm{Pic}^0(X)
\]
is surjective
since $\varphi_{A_1^{\otimes m}}$ is a surjective map and is the composition of $\varphi_{A_1}$ and an \'{e}tale endomorphism
\[
\mathrm{Pic}^0(X)\ni [M]\mapsto [M^{\otimes m}]\in\mathrm{Pic}^0(X).
\]
This is equivalent to the existence of an isomorphism $\xi\in\mathrm{Aut}_0(X,\Delta)$ such that $\xi^*A_1\sim A_2$ for any $A_2$ algebraically equivalent to $ A_1$.
We complete the proof.
\end{proof}

To prove Theorem \ref{thm--proj-our-moduli}, we show the following by applying Lemma \ref{lem-pic-aut}.

\begin{prop}\label{prop--pic-aut-unif-ad}
    Let $f\colon(X,\Delta,A)\to \mathbb{P}^1$ be a uniformly adiabatically {\rm K}-stable klt-trivial fibration with $-(K_X+\Delta)$ not numerically trivial but nef.

    Then $\mathrm{dim}\,\mathrm{Aut}_0(X,\Delta)=\mathrm{dim}\,\mathrm{Pic}^0(X)$ and for any two ample line bundles $A_1$ and $A_2$ algebraically equivalent to each other, there exists $\varphi\in\mathrm{Aut}_0(X,\Delta)$ such that $\varphi^*A_1\sim A_2$.
\end{prop}

\begin{proof}
    Take $s\in\mathbb{Q}_{>0}$ such that $-(K_X+\Delta)\sim_{\mathbb{Q}}sf^*\mathcal{O}(1)$.
Then we see by \cite[Theorem 1.1]{Hat} that for any distinct three closed points $p_1,p_2,p_3\in\mathbb{P}^1$, $(X,\Delta+\frac{s}{3}\sum_{i=1}^3f^{-1}(p_i))$ is klt and $K_X+\Delta+\frac{s}{3}\sum_{i=1}^3f^{-1}(p_i)\sim_{\mathbb{Q}}0$.
We claim that 
\begin{equation}\label{eq--prop--pic-aut}
\mathrm{Aut}_0(X,\Delta)=\mathrm{Aut}_0\left(X,\Delta+\frac{s}{3}\sum_{i=1}^3f^{-1}(p_i)\right).
\end{equation}
Indeed, $\mathrm{Aut}_0(X,\Delta)$ acts on $\mathbb{P}^1$ but $\mathrm{Aut}_0(X,\Delta)$ is an Abelian variety by \cite{CM}.
Let $G$ be the image of the group homomorphism $\mathrm{Aut}_0(X,\Delta)\to PGL(2)$.
Since $G$ is a proper linear algebraic group, $G$ is a point. 
Thus, $\mathrm{Aut}_0(X,\Delta)$ fixes $f^{-1}(p)$ for any $p\in\mathbb{P}^1$ and (\ref{eq--prop--pic-aut}) holds.
By Lemma \ref{lem-pic-aut}, 
\[
\mathrm{dim}\,\mathrm{Aut}_0\left(X,\Delta+\frac{s}{3}\sum_{i=1}^3f^{-1}(p_i)\right)=\mathrm{dim}\,\mathrm{Pic}^0(X).
\]
Thus, we complete the proof of the first assertion by (\ref{eq--prop--pic-aut}).
The second assertion follows in the same way as Lemma \ref{lem-pic-aut}.
\end{proof}

\begin{proof}[Proof of Theorem \ref{thm--proj-our-moduli}]
By \cite[Proposition 2.7]{kollar-moduli-stable-surface-proj}, there exist a proper normal variety $B'$, a finite surjective morphism $g\colon B'\to B$ and a morphism of stacks $\tilde{g}\colon B'\to\mathscr{M}_{d,v,u,r}$ such that $\pi\circ\tilde{g}=\iota\circ g$, where $\iota\colon B\to M_{d,v,u,r}$ and $\pi\colon \mathscr{M}_{d,v,u,r}\to M_{d,v,u,r}$ are the natural morphisms.
We set $w$ as in Theorem \ref{quesmain}.
Let $f\colon(X,A)\to B'$ be the pullback of the universal family $(\mathscr{U},\mathscr{A})$ on $\mathscr{M}_{d,v,u,r}$ via $\tilde{g}$ (cf.~\cite[Remark 6.5]{HH}) with $\mathrm{vol}(A_{b'})=w$ for any point $b'\in B'$.
Then, $A$ is $f$-ample and $(X_b,A_b)$ is specially K-stable for any closed point $b\in B'$.

Here, we claim that the family $X\to B'$ has maximal variation.
To show this, assume the contrary and that there exists a proper curve $C\subset B'$ such that $C$ passes through a very general point and for any general two closed points $p_1,p_2\in C$, $X_{p_1}$ and $X_{p_2}$ are isomorphic.
Then, $A_p$ and $A_q$ are algebraically equivalent for any very general two closed points $p,q\in C$.
By Proposition \ref{prop--pic-aut-unif-ad}, we see that $C$ is contained in a fiber of $g$.
This contradicts to the finiteness of $g$.
Thus, the family $X\to B'$ has maximal variation.

The CM line bundle $\lambda_{\mathrm{CM},f}=g^*(\Lambda_{\text{CM},w}|_B)$ on $B'$ is big and nef by Theorem \ref{thm--main--2}.
Thus, we have that $(\Lambda_{\mathrm{CM},w}|_B)^{\mathrm{dim}\,B}>0$.
By the Nakai-Moishezon criterion \cite[Theorem 3.11]{kollar-moduli-stable-surface-proj}, $\Lambda_{\text{CM},w}|_B$ is ample and hence $B$ is projective.
\end{proof}


\end{document}